\documentclass[11pt]{amsart}
\usepackage{amsmath,amssymb}

\usepackage{fig4tex}

\newtheorem{thm}{Theorem}[section]

\newtheorem{lem}[thm]{Lemma}
\newtheorem{prop}[thm]{Proposition}

\newtheorem{hyp}[thm]{Assumption}
\newtheorem{notn}[thm]{Notation}

\theoremstyle{definition}

\theoremstyle{remark}
\newtheorem{rem}[thm]{Remark}

\newcommand{\iti}[1]{\noindent\mbox{{\em (\romannumeral #1)\/}}}
\newcommand{\db}[1]{_{\raise-0.3ex\hbox{$\scriptstyle #1$}}}
\newcommand{\dd}[1]{_{\raise-1.5pt\hbox{$\scriptstyle #1$}}}

\newcommand{\di}{\displaystyle}

\newcommand{\dr}{{\rm d}}

\newcommand  {\N}{{\mathbb N}}

\newcommand  {\R}{{\mathbb R}}

\newcommand {\Id}{\mathbb {I}}

\renewcommand  {\H}{{\mathrm H}}

\newcommand{\BB}{\boldsymbol{\mathsf B}}

\newcommand  {\LL}{\boldsymbol{\mathsf L}}
\renewcommand  {\L}{{\mathrm L}}

\newcommand  {\TT}{\boldsymbol{\mathsf T}}

\renewcommand  {\dd}{{\boldsymbol{\mathsf d}}}

\newcommand  {\eps}{\varepsilon}
\newcommand  {\ff}{{\boldsymbol{\mathsf f}}}

\newcommand  {\nn}{\boldsymbol{\mathsf n}}
\newcommand  {\nnu}{\boldsymbol{\mathsf \nu}}

\newcommand  {\s}{\mathsf s}
\newcommand  {\ttt}{\boldsymbol{\mathsf t}}
\newcommand  {\uu}{\boldsymbol{\mathsf u}}
\newcommand  {\vv}{\boldsymbol{\mathsf v}}

\newcommand  {\xx}{\boldsymbol{\mathsf x}}

\newcommand  {\cH}{\mathcal{H}}

\newcommand  {\cJ}{\mathcal{J}}
\newcommand  {\cK}{\mathcal{K}}

\newcommand  {\bH}{\mathbf{H}}
\newcommand  {\bL}{\mathbf{L}}
\newcommand  {\bV}{\mathbf{V}}

\newcommand {\su}{\mathsf{u}}

\renewcommand {\sp}{\mathsf{p}}
\newcommand {\sq}{\mathsf{q}}
\renewcommand {\sf}{\mathsf{f}}
\newcommand {\sg}{\mathsf{g}}

\newcommand  {\f}{\mathfrak{f}}
\newcommand  {\fg}{\mathfrak{g}}
\newcommand  {\fp}{\mathfrak{p}}
\newcommand  {\fq}{\mathfrak{q}}

\newcommand{\B}{\mathsf B}

\newcommand{\sN}{\mathsf N}

\newcommand  {\F}{\mathfrak{F}}
\begin{document}
%%-----------------------------
%%      the top matter
%%-----------------------------
\title[Equivalent Conditions for Elasto-Acoustics]{Equivalent Boundary Conditions for an Elasto-Acoustic Problem set in a Domain with a Thin Layer.}
%\thanks{INRIA Bordeaux Sud-Ouest}
\thanks{This work was supported by the project HPC GA 295217 - IRSES2011}% At most 5 thanks
\author{Victor P\'eron}\address{ LMAP CNRS UMR 5142 \& Team MAGIQUE 3D INRIA, Universit\'e de Pau et des Pays de l'Adour, France}
%\author{...}\address{...}
%\author{...}\address{...}
%
%\date{...}
%
\begin{abstract}
We present equivalent conditions and asymptotic models for the diffraction problem of elastic and acoustic waves in a solid medium surrounded by a thin layer of fluid medium. Due to the thinness of the layer with respect to the wavelength, this problem is well suited for the notion of equivalent conditions and the effect of the fluid medium on the solid is as a first approximation local. We derive and validate equivalent conditions up to the  fourth order for the elastic displacement. These conditions approximate the acoustic waves which propagate in the fluid region. This approach leads  to solve only elastic equations. The construction of equivalent conditions is based on a multiscale expansion in power series of the thickness of the layer for the solution of the transmission problem. 
\end{abstract}
%
%\begin{resume}
%On pr\'esente des conditions d'imp\'edances ainsi que des mod\`eles asymptotiques adapt\'es \`a un probl\`eme de transmission pour un couplage \'elasto-acoustique dans un milieu avec une couche mince acoustique. Ce type de probl\`eme est bien pos\'e pour la notion de conditions d'imp\'edances car la faible \'epaisseur de la couche assure que l'effet acoustique sur la  d\'eformation \'elastique est en premi\`ere approximation locale. On justifie rigoureusement ces conditions d'imp\'edances. Cette m\'ethode permet d'approcher la solution du probl\`eme exact par celle d'un probl\`eme \'elastique avec une condition d'imp\'edance. L'erreur d'approximation est contr\^ol\'ee par rapport \`a l'\'epaisseur de la couche.
%\end{resume}
%
\maketitle
%%-----------------------------
%%      your text
%%-----------------------------
\section*{Introduction}

  	The concept of Equivalent Boundary Conditions (also called approximate, effective, or impedance conditions) is classical in the 
 modeling of wave propagation phenomena. Equivalent Conditions (ECs) are usually introduced to reduce the computational domain. The main idea consists to replace an ``exact'' model inside a part of the domain (for instance a thin layer of dielectric %
  or a highly absorbing material) by an approximate boundary condition. This idea is pertinent when the effective condition can be readily handled for numerical computations, for instance when this condition is local~\cite{EN93,SeVo95,BL96,HJN08}. In the 1990's Engquist--N{\'e}d{\'e}lec~\cite{EN93}, Abboud--Ammari~\cite {AA95}, Bendali--Lemrabet~\cite {BL96}, Ammari--N{\'e}d{\'e}lec~\cite {AN98},  and Lafitte~\cite {L99} derived equivalent conditions for acoustic and electromagnetic scattering problems to approximate an obstacle coated by a thin layer of dielectric (absorbing) material inside the domain of interest. Impedance conditions are also used to reduce the computational domain for scattering problems from highly absorbing obstacles~\cite{Rytov40,Leo48,LL93,SeVo95,ABV01,HJN08}.

	 The main application of this work concern the mathematical modeling of earthquake on the Earth's surface. %This work enters into the scope of projects which involve high performance computing for geophysical applications. 
The simulation of large-scale geophysics phenomena represents a main challenge for our society. Seismic activities worldwide have shown how crucial it is to enhance our understanding of the impact of earthquakes. In this context, the coupling of elastic and acoustic waves equations is essential if we want to reproduce real physical phenomena such as an earthquake. %Acoustic waves propagate into the ocean whereas elastic waves propagate into the Earth subsurfaces. 
We can thus take into account the effects of the ocean on the propagation of seismic waves. Elasto-acoustic coupling problems are rather classical in the mathematical modeling of wave propagation phenomena. We refer to several works in Ref.~\cite {J83,LM95,DN00,NSW00,HKM08,MS09,Am10} and the monography~\cite[\S 5.4.e]{CDN10} which concern the direct problem for fluid-structure interaction systems and acoustic (or elastic) scattering by smooth elastic obstacles. %We refer also to~\cite{MS09} where the authors study in particular a direct problem for fluid-solid interactions. 
We intend to work in the context of this application for which we consider that the medium consists of land areas surrounded by fluid zones whose thickness is very small, typically with respect to the wavelength. This raises the difficulty of applying a finite element method on a mesh that combines fine cells in the fluid zone and much larger cells in the solid zone. To overcome this difficulty and to solve this problem, we adopt an asymptotic method which consists to ``approximate'' the fluid portion by an equivalent boundary condition. This boundary condition is then coupled with the elastic wave equation and a finite element method  can be applied to solve  the resulting boundary value problem. 
		
	In this paper, we present elements of derivation together with mathematical justifications for equivalent boundary conditions, which appear as a first, second, third or fourth order approximations (Sec.~\ref{Ecs}) with respect to the small parameter (the thickness of the fluid layer) and satisfied by the elastic displacement $\uu$. This work is concerned essentially with theoretical objectives. The numerical pertinence of these ECs up to the third order have already been shown for the two-dimensional problem~\cite{DP13}.	
	
	There are several similarities in this paper and in the works in Ref.~\cite{EN93,AA95,BL96} in which the authors considered the problem of a time-harmonic wave for the Helmholtz equation. In Ref.~\cite{EN93},  the authors construct ECs up to order $3$ whereas in Ref.~\cite{BL96} the authors derive and analyze ECs up to order $4$. In this paper, we construct and analyze ECs up to order $4$ for elasto-acoustics. As in Ref.~\cite{EN93,BL96}, the construction of the ECs relies on a multiscale expansion of the exact solution in power series of the small parameter. %The multiscale analysis to describe a transmission problem with a thin layer is rather classical~\cite{CCDV06,DPP11,ST11}. 
	However, there are several differences between the results of this paper and the works in Ref.~\cite{EN93,BL96} since ECs are not of the same nature%in this work  and in~\cite{EN93,BL96}
	. The ECs are of ``$\uu\cdot\nn$--to--$\TT(\uu)$'' nature for elasto-acoustics since a local impedance operator links the normal traces of $\uu$ and the stress vector $\TT(\uu)$ whereas the ECs are of ``Dirichlet--to--Neumann'' nature for acoustics in Ref.~\cite{EN93,BL96}. We compare impedance operators for elasto-acoustics and acoustics~\cite{BL96} precisely. One difficulty to validate the equivalent conditions lies in the proof of uniform energy estimates for the exact solution of the elasto-acoustic coupling. We overcome this difficulty using a compactness argument and removing a discret set of resonant frequencies (which may appear in the solid part of the domain). We revisit the proof of uniform estimates in Ref.~\cite{BL96}. We can not apply straightforwardly the proof of estimates in Ref.~\cite{BL96} to the elasto-acoustic coupling since the transmission conditions (which are purely natural) for elasto-acoustics appear in the weak formulation and play a crucial role. We prove well-posedness and convergence results for ECs up to the fourth order. 
		
	The outline of the paper proceeds as follows. In Section \ref{SMM}, we introduce the mathematical model and the framework for the elasto-acoustic problem and equivalent conditions. We present briefly a formal derivation of equivalent conditions. 	Then we present uniform estimates for the solution of the transmission problem. In Section \ref{Smr}, we present equivalent conditions and asymptotic models associated with the solution of the exact problem. We eventually compare these equivalent conditions with the works~\cite{BL96}. In Section \ref{Sue}, we prove uniform estimates for the solution of the elasto-acoustic problem. In Section \ref{SDEC}, we derive and validate a two-scale asymptotic expansion at any order for the solution of the problem, and we construct formally ECs. In Section \ref{SAEC}, we prove stability results for ECs and the convergence of ECs towards the exact model.

\section{The Mathematical Model}
\label{SMM}
 In this section, we introduce the model problem (\S \ref{Smp}) and the framework for the elasto-acoustic problem. Then we remind  the definition of equivalent conditions and we state uniform estimates for the solution of the exact problem. We start this section with a formal derivation of the approximate boundary conditions. 
 
 \subsection{Formal derivation of equivalent conditions}

In this section, we present briefly a formal derivation of equivalent conditions. We summarize this process in two steps.  All the details and formal calculi are presented in Sec. \ref{SDEC}.
 
 \subsubsection*{First step : a multiscale expansion}
The first step consists to derive a multiscale expansion for the solution $(\uu_\eps,\sp_{\eps})$ of the model problem \eqref{EA} (Sec. \ref{Smp}): it possesses  an asymptotic expansion in power series of the small parameter $\eps$
\begin{gather*}
%\label{Euu}
   \uu_{\eps}(\xx) = \uu_0(\xx) + \eps\uu_1(\xx) + \eps^2 \uu_2(\xx)+ \cdots  \quad\mbox{in}\quad \Omega_{\mathsf{s}}
   %\mathcal{O}(\eps^3)
   \, ,\\[0.8ex]
   \label{Esp}
   \sp_{\eps}(\xx) = \sp_0(\xx;\eps) + \eps \sp_1(\xx;\eps) + \eps^2 \sp_2(\xx;\eps)  +\cdots  \quad \mbox{in}\quad \Omega^{\eps}_\mathsf{f}
   %+ \mathcal{O}(\eps^2)
     \, ,\\
   \quad   \mbox{with}\quad
   \sp_j(\xx;\eps) = \fp_j(y_{\alpha},\frac{y_{3}}{\eps})\, .
\end{gather*}
Here $\xx\in\R^3$ are the cartesian coordinates and $ (y_{\alpha}, y_{3})$ is a ``{\em normal coordinate system}'' ~\cite {CDFP11,Pe13} to the surface $\Gamma=\partial\Omega_{\mathsf{s}} $ on the manifold $\Omega^{\eps}_\mathsf{f}$ : $y_{\alpha}$ ($\alpha\in\{1,2\}$) is a tangential coordinate on $\Gamma$ and $y_{3} \in(0,\eps)$ is the distance to the surface $\Gamma$.
The term $\fp_j$ is a ``{\em profile}'' defined on $\Gamma\times (0,1)$. Formal calculi  are presented in Sec. \ref{SFAE} and the first terms $(\fp_j,\uu_j)$ for $j=0,1,2,3$ are explicited in Sec. \ref{ft}.

 \subsubsection*{Second step : construction of equivalent conditions}

The second step consists to identify for all $k\in\{0,1,2,3\}$ a simpler problem satisfied by the truncated expansion 
$$
\uu_{k,\eps}:=\uu_0 + \eps\uu_1 + \eps^2 \uu_2+ \cdots+ \eps^k \uu_k 
$$
 up to a residual term in $\mathcal{O}(\eps^{k+1})$. The simpler problem writes
\begin{equation*}
 \left\{
   \begin{array}{lll}
 \nabla \cdot  \underline{\underline{\sigma}} (\uu_\eps^{k})+\omega^2 \rho \uu_\eps^{k} =\ff
 \quad&\mbox{in}\quad \Omega_{\mathsf{s}}
\\[0.8ex]
\TT(\uu_\eps^{k})+\B_{k,\eps}(\uu_\eps^{k}\cdot\nn)\, \nn =0% \quad \mbox{and} \quad \uu \cdot\nn=0
  \quad &\mbox{on} \quad \Gamma\ .
 %\ 
   \end{array}
    \right.
\end{equation*}
Here $\ff$  is the data of the model problem \eqref{EA} and $\B_{k,\eps}$ is a surfacic differential operator acting on functions defined on $\Gamma$ and which depends on $\eps$
\begin{gather*}
\label{B0}
\B_{0,\eps}=0
   \, ,\\[0.8ex]
   \label{B1}
   \B_{1,\eps}(\su)= -\eps \omega^2 \rho_{\sf} \, \su
  \, ,\\[0.8ex]
   \label{B2}
\B_{2,\eps}(\su)= - \eps   \omega^2 \rho_{\sf} \left( 1 -\eps \cH(y_\alpha) \right)\su
\, ,
\\[0.8ex]
\B_{3,\eps}(\su)=-\eps   \omega^2 \rho_{\sf}\left(1 -\eps \cH(y_\alpha) +\frac{\eps^2}{3} \left[ \Delta_{\Gamma} + \kappa^2 \Id  + 4 \cH^2 (y_\alpha) -\cK(y_\alpha) \right]  \right)\su 
\, .
\end{gather*}
Here $\cH$ and  $\cK$ denote the \textit{mean curvature} and the \textit{Gaussian  curvature} of the surface $\Gamma$  and $\Delta_{\Gamma}$ is the Laplace-Beltrami operator along $\Gamma$. Equivalent conditions are stated in Sec.  \ref{S3DECs}. 
The construction of these conditions is detailed in Sec. \ref{CECs}.

\subsection{The model problem}
\label{Smp}
Our interest lies in an elasto-acoustic wave propagation problem  in time-harmonic regime set in a domain with a thin layer. We consider the fluid-solid transmission problem
\begin{equation}
\label{EA}
 \left\{
   \begin{array}{lll}
  \Delta \sp_{\eps}+\kappa^2 \sp_{\eps} = 0   \quad&\mbox{in}\quad \Omega^{\eps}_\mathsf{f}
\\[0.5ex]
 \nabla \cdot   \underline{\underline{\sigma}}( \uu_{\eps})+\omega^2 \rho \uu_{\eps} =\ff
 \quad&\mbox{in}\quad \Omega_{\mathsf{s}}
\\[0.5ex]
\partial_{\nn }\sp_{\eps}=\rho_{\sf} \omega^2 \uu_{\eps}\cdot \nn \quad&\mbox{on}\quad \Gamma
\\[0.5ex]
\TT(\uu_{\eps}) = -\sp_{\eps} \nn  \quad &\mbox{on} \quad \Gamma
\\[0.2ex]
%\mbox{external B.C.} 
\sp_{\eps}=0 
  \quad &\mbox{on}\quad \Gamma^\eps \ , % \partial_{\nn }\sp_{\eps} - i \kappa\sp_{\eps}=0 
   \end{array}
    \right.
\end{equation}
set in a smooth bounded simply connected domain $\Omega^\eps$ in $\R^3$ made of a solid, elastic object occupying  a smooth connected subdomain $\Omega_{\mathsf{s}}$ entirely immersed in a fluid region occupying the subdomain $\Omega^{\eps}_\mathsf{f}$. The domain $\Omega^{\eps}_\mathsf{f}$ is a thin layer of uniform thickness $\eps$ (i.e. the euclidean distance between surfaces $\Gamma$ and $\Gamma^\eps$ is $\eps$),  see figure \ref{FCylGeo}. We denote by $\Gamma^\eps $ the boundary of the domain $\Omega^\eps$, and by $\Gamma$ the interface between the subdomains  $\Omega^{\eps}_\mathsf{f}$ and $ \Omega_{\mathsf{s}}$. We denote by $\nn$ the unit normal to $\Gamma$ oriented from $\Omega_{\mathsf{s}}$ to  $\Omega^{\eps}_\mathsf{f}$.
\input contribF4T.tex
\begin{figure}[ht]
%%%%%%%%%%%%%%%%%%%%%%%%%%%%%%%%
%%%%%%%% DEBUT DESSIN
%%%%%%%%%%%%%%%%%%%%%%%%%%%%%%%
\def\EpsTube{11}
\begin{center}
\psset(width=2)
\figinit{1pt}
\pssetupdate{yes}
% pt de la courbe \Sigma?
\figpt 8:(-48,-18)%\figpt 9:(-20,50)\figpt 10:(50,20)
%\figpt 12:(10,10)\figpt 13:(0,-10)
%\figpt 14:(-20,15)%
% pt de la courbe \Sigma_{h}
\figpt 22:(-92,25.5)\figpt 23:(-27.2,60.8)
%\figpt 23:(-37.2,30.8)
\figpt 24:(46.3,41.5)\figpt 25:(70.8,2.8)
\figpt 26:(29,-31)
%\figpt 27:(-27.2,8)
\figpt 27:(-20,1)\figpt 28:(-77,-11)
% pt de la courbe \Sigma
\figpt 1:(-28.3,71.9)
\figpt 2:(-12,-4.5)
\figpt 3:(50,50)
%\nn
\figpt 4:(58,62)\figpt 5:(30,-40)
%(46.3,41.5)
%\figpt 6:(-30,0)\figpt 7:(-80,-20)
% ecriture de \Sigma, \Omega\iso, \Omega\con et \Partial\Omega, et \nn et h et ${\mathcal{O}}$
\figpt 16:(-70,20)\figpt 15:(0,30)
\figpt 17:(40,10) \figpt 18:(-60,-20)
\figpt 21:(50,65) \figpt 31:(30,-18)  \figpt 34:(-32,65) 
 \figpt 35:(-89,20)
   \figpt 36:(-103,5)
% pt vecteurs normaux
\figpt 19:(61,65) \figpt 20:(3,27)  \figpt 30:(42.6,33)  \figpt 32:(55,7) \figpt 33:(61.6,6.1) 
%%%%%%%%%%%%%%%%%%%%%
%% fichier graphique
%%%%%%%%%%%%%%%%%%%
\psbeginfig{}
\psset arrowhead(fillmode=no,length=3)\psarrow[1,23]
\psset arrowhead(fillmode=no,length=3)\psarrow[23,1]
%\psarrow[23,2]
%\pssetfillmode{yes}\pssetgray{0.5}
%\pscurve[1,2,3,4,5,6,7,1,2,3]
%\psBezier 4[1,2,3,4,5,6,7,1]
%\pscurve[8,9,10,12,13,14,8,9,10]
\pssetfillmode{no}\pssetgray{0}
%\pscurve[8,9,10,12,13,14,8,9,10]
%\pssetfillmode{yes}\pssetgray{0.5}
\pscurve[22,23,24,25,26,27,28,22,23,24]
\figptscontrolcurve 40,\NbC[22,23,24,25,26,27,28,22,23,24]
%\psreset{first}\psset (width=1.9)
\psEpsLayer \EpsTube,\NbC[40,41,42,43,44,45,46,47,48,49,50,51,52,53,54,55,56,57,58,59,60,61]
%\pscurve[40,43,46,49,52,55,58,40,43]
%\figptscontrolcurve 70,\NbC[40,43,46,49,52,55,58,40,43]
%\psreset{first}\psset (width=1.9)
%\psEpsLayer \EpsTubebis,\NbC[70,71,72,73,74,75,76,77,78,79,80,81,82,83,84,85,86,87,88,89,90,91]
%\pscurve[22,23,24,25,26,27,28,22,23,24]
%\figptscontrolcurve 70,\NbC[22,23,24,25,26,27,28,22,23,24]
%\psreset{first}\psset (width=1.9)
%\psEpsLayer \EpsTubebis,\NbC[70,71,72,73,74,75,76,77,78,79,80,81,82,83,84,85,86,87,88,89,90,91]
%\'epaisseur h
%\psline[23,2]
%\psarrow[23,2]
%\psarrow[2,23]
%\pssetfillmode{yes}\pssetgray{0.5}
%Vecteurs
%\psarrow[32,25]
\psset arrowhead(fillmode=yes,length=3)\psarrow[24,4]
\pssetfillmode{yes}\pssetgray{0.5}
\psendfig
%%%%%%%%%%%%%%%%%%%
%% writing
%%%%%%%%%%%%%%%%%%%
\figvisu{\figBoxA}{
%{Figure 1}\ --\ Un voisinage tubulaire $\mathcal{O^{\prime}}\subset\mathcal{O}$ de $\Sigma$ dans $\Omega$
}
{
\figwritew 15: $\Omega_\s$ %${\Omega\con}$
(6pt)
%\figwritec [16]{$\Omega\iso$}
\figwritew 2: $\Omega^{\eps}_\mathsf{f}$ (4pt)
\figwritec [34]{$\eps$}
\figwritew 36: $\Gamma^\eps$ (2pt)
%\figwriten 31: $\Sigma_{h}$ (2pt)
\figwritew 4: $\nn$(4pt)
%\figwritew 32: $\nn$(3pt)
%\figwritee 20: $\nn$(3pt)
\figwrite [35]{$\Gamma$}
\figsetmark{$\figBullet$}
%\figwritep[1,2,3,4,5,6,7,22,23,24,25,26,27,28]
}
\centerline{\box\figBoxA}
 \caption{A cross-section of the domain $\Omega^\eps$ and its subdomains $\Omega_\mathsf{s}$ and $\Omega^{\eps}_\mathsf{f}$ }
\label{FCylGeo}
\end{center}
\end{figure}

In the elasto-acoustic system \eqref{EA}, we denote the unknowns by $\uu_{\eps}$ for the elastic displacement and by $\sp_{\eps}$ for the acoustic pressure. The time-harmonic wave field with angular frequency $\omega$ is characterized by using the Helmholtz equation for the pressure $\sp_{\eps}$, and by using an anisotropic discontinuous linear elasticity system for the displacement $\uu_{\eps}$%in time-harmonic regime
. These equations contain several physical constants: $\kappa=\omega/c$ is the acoustic wave number, $c$ is the speed of the sound, $\rho$ is the density of the solid, and $\rho_{\sf}$ is the density of the fluid.  All these constants are independent of $\eps$. 

 In the linear elastic equation, $\nabla \cdot$ denotes the divergence operator for tensors, $ \underline{\underline{\sigma}}( \uu)$  is the stress tensor given by Hooke's law 
$$
\underline{\underline{\sigma}}( \uu)=\underline{\underline{C}}\ \underline{\underline{\epsilon}} (\uu)\ .
$$
 Here $ \underline{\underline{\epsilon}}  (\uu )=\frac12(\underline{\underline{\nabla}} \uu+\underline{\underline{\nabla}} \uu^T )$ is the strain tensor where $\underline{\underline{\nabla}}$ denotes the gradient operator for tensors, and $\underline{\underline{C}} =\underline{\underline{C}} (\xx)$ is the elasticity tensor. The components of $\underline{\underline{C}} $ are the elasticity moduli $C_{ijkl}$ : $\underline{\underline{C}} =\left(C_{ijkl}(\xx)\right)$. The \emph{traction operator}  $\TT$ is a surfacic differential operator defined on $\Gamma$ as 
 $$
 \TT(\uu)=
%\underline{\underline{C}} \ \underline{\underline{\epsilon}} 
 %\underline{\underline{C}} \ \underline{\underline{\nabla}} \uu  \,\nn 
 \underline{\underline{\sigma}}(\uu)\nn \ .
$$ 
The right-hand side $\ff$ is a data with support in $\Omega_{\s}$. The first transmission condition set on $\Gamma$ is a kinematic interface condition whereas the second one is a dynamic interface condition. The kinematic condition requires that the normal velocity of the fluid match the normal velocity of the solid  on the interface $\Gamma$. The dynamic condition results from the equilibrium of forces on the interface $\Gamma$. The transmission conditions are natural.

\begin{rem}
We consider in this paper mainly Dirichlet external boundary conditions. In~\cite[Appendix A]{Pe13}, we present also equivalent conditions up to the second order for the elato-acoustic problem complemented with a  Fourier-Robin boundary condition set on $\Gamma^\eps$.
% $$
 %\quad \partial_{\nn }\sp_{\eps} - i \kappa \sp_{\eps}=0 \quad \mbox{on} \quad \Gamma^\eps\ .
%$$
\end{rem}

\medskip
In the framework above we address the issue of Equivalent Conditions (ECs) for the elastic displacement $\uu_{\eps}$  as $\eps \to 0$, see Section \ref{Ecs}.
%\begin{enumerate}
%\item The issue of Equivalent Conditions (ECs) for the elastic displacement $\uu_{\eps}$  as $\eps \to 0$, see Section \ref{Ecs},
%\item The issue of Uniform Estimates for the displacement $\uu_{\eps}$ and the pressure $\sp_{\eps}$ solutions of the problem \eqref{EA} as $\eps \to 0$, see Section \ref{UEs}.
%\end{enumerate} 
This issue is linked with the question of Uniform Estimates for the couple $(\uu_{\eps},\sp_{\eps})$ solution of the problem \eqref{EA} as $\eps \to 0$ (Section \ref{UEs}) since it is a main ingredient in the mathematical justification of ECs. To answer these questions, we make hereafter several assumptions on the data and on the regularity of the surface $\Gamma$. %These assumptions  simplify the asymptotic modeling.

\subsection{Framework}
\label{Sa}

We will work under usual assumptions (symmetry and positiveness) on the elasticity tensor. %$ \underline{\underline{C}} $ .
\begin{hyp}
\label{AC}
\iti1 The elasticity moduli $C_{ijkl}(\xx)$ are real valued smooth functions in $\overline{\Omega_{\s}}$.

\iti2 The tensor $\underline{\underline{C}} $ is symmetric : 
$$
  C_{ijkl}=C_{jikl}=C_{klij} \quad\mbox{almost everywhere in} \, \Omega_{\s} \ .
$$
\iti3  The tensor $\underline{\underline{C}} $ is  positive : 
$$\exists \alpha >0  , \quad \forall \xi=(\xi_{ij}) \,\mbox{symmetric tensor}, \,\sum_{i,j,k,l} C_{ijkl} \xi_{ij} \overline{\xi_{kl}} \geqslant \alpha \sum_{i,j} |\xi_{ij}|^2\ .
$$
 \end{hyp}

\begin{rem} 
According to the assumption \ref{AC} \iti2, the Hooke's law writes also $\underline{\underline{\sigma}}( \uu)=\underline{\underline{C}}\ \underline{\underline{\nabla}} \uu$ . Hence, the assumption \ref{AC} \iti3 ensures that the matrix differential operator $\nabla \cdot  \underline{\underline{\sigma}}+\omega^2 \rho \Id $
is strongly elliptic.
\end{rem}

Some resonant frequencies  may appear in the solid domain. However, we prove uniform estimates  for the elasto-acoustic field $(\uu_{\eps},\sp_{\eps})$ as well as ECs for $\uu_{\eps}$ when $\eps\to 0$ under the following  spectral assumption on the limit problem set in the solid part $\Omega_{\mathsf{s}}$, and when $\ff=0$. 
\begin{hyp}
\label{Hw}
 The angular frequency $\omega$ is not an eigenfrequency of the problem 
 \begin{equation*}
\label{SA}
 \left\{
   \begin{array}{lll}
\nabla \cdot  \underline{\underline{\sigma}}(\uu)+\omega^2 \rho \uu =0
 \quad&\mbox{in}\quad \Omega_{\mathsf{s}}
\\[0.5ex]
\TT(\uu) =0 \quad &\mbox{on} \quad \Gamma \ .
   \end{array}
    \right.
\end{equation*}
 \end{hyp}

 Our whole analysis is valid under the following assumption on the surfaces $\Gamma$ and $\Gamma_{\eps}$.
\begin{hyp}
\label{HG}
The fluid-solid interface $\Gamma$ and the surface $\Gamma_{\eps}$ are smooth.
 \end{hyp}

For the sake of simplicity in the asymptotic modeling, we will work under the following assumption on the data $\ff$.
\begin{hyp}%[Regularity on the data $\sp_{i}$]
\label{Hspi}
The right-hand side $\ff$ in \eqref{EA} is a smooth $\eps$-independent data. %with support in $\Omega_{\s}$.
%, $\sp_{i}$ is supposed to be a "smooth" data.
% For the sake of simplicity in the asymptotic modeling, we assume that the fluid-solid interface $\Gamma$ is smooth.
 \end{hyp}

In the framework above, we prove in this paper that it is possible to replace 
%our aim consists in replacing
 the fluid region $\Omega^{\eps}_\mathsf{f}$ 
 by appropriate boundary conditions called equivalent conditions and  set on $\Gamma$. 
 %This is the main result of the paper  hereafter. %We assume that the small parameter is $\eps\ll\lambda$.

\subsection{Validation of equivalent conditions}
\label{Ecs}
In this paper, %we derive and prove equivalent conditions (ECs) :
 we derive surfacic differential operators $\B_{\eps}$
$$
\B_{\eps} :\mathcal{C}^\infty(\Gamma) \rightarrow  \mathcal{C}^\infty(\Gamma) \ , 
%\Rd \quad \bH^1(\Omega_{\mathsf{s}}) \rightarrow  \bH^\frac12(\Gamma)\Bk
$$
together with $\tilde\uu_\eps$ which is a solution of the boundary value problem 
\begin{equation}
\label{Eec}
 \left\{
   \begin{array}{lll}
 \nabla \cdot  \underline{\underline{\sigma}} (\tilde\uu_\eps)+\omega^2 \rho \tilde\uu_\eps =\ff
 \quad&\mbox{in}\quad \Omega_{\mathsf{s}}
\\[0.5ex]
\TT(\tilde\uu_\eps)+\B_{\eps}(\tilde\uu_\eps\cdot\nn) \nn = 0 
  \quad &\mbox{on} \quad \Gamma\ .
 %\ 
   \end{array}
    \right.
\end{equation}
 Then in the framework of Sec. \ref{Sa}, we prove uniform estimates for the error between the exact solution $\uu_{\eps}$ in \eqref{EA} and $\tilde\uu_\eps$ provided $\eps$ is small enough: 
\begin{equation}
\label{ek}
\|  \uu_{\eps}- \tilde\uu_\eps \|_{1,\Omega_{\mathsf{s}}} \leqslant C \eps^{k+1}\ , 
\end{equation}
with $k\in\{0,1,2,3\}$, see Th. \ref{CVec} for the main result and precise estimates. Here, we denote by $\| \cdot \|_{1,\Omega_{\mathsf{s}}}$  the norm in the Sobolev space $\bH^1(\Omega_{\mathsf{s}})=\H^1(\Omega_{\mathsf{s}})^3$. We say that the equivalent condition is of order $k+1$ when such an a priori estimate \eqref{ek} holds. Then we define $\uu_\eps^{k}=\tilde\uu_\eps$  and we denote by $\B_{k,\eps}$ the operator $\B_{\eps}$ corresponding to the order $k+1$, Sec. \ref{Smr}. The validation of ECs relies on uniform estimates for solutions $(\uu_{\eps},\sp_{\eps})$ of \eqref{EA} as $\eps\to 0$. This issue is developped in Sec. \ref{UEs}.

\subsection{Uniform estimates}
\label{UEs}
We introduce a suitable variational framework for the solution of the problem \eqref{EA} with more general right-hand sides. This framework is useful to prove error estimates \eqref{ek}.%$(\sf,\sg,\hh)$.

\subsubsection*{Weak solutions}
For given data $(\ff,\sf,\sg)$ we consider the boundary value problem
\begin{equation}
\label{EAfgh}
 \left\{
   \begin{array}{lll}
  \Delta \sp_{\eps}+\kappa^2 \sp_{\eps} = \sf  \quad&\mbox{in}\quad \Omega^{\eps}_\mathsf{f}
\\[0.5ex]
 \nabla \cdot  \underline{\underline{\sigma}}(\uu_{\eps})+\omega^2 \rho \uu_{\eps} =\ff
 \quad&\mbox{in}\quad \Omega_{\mathsf{s}}
\\[0.5ex]
\partial_{\nn }\sp_{\eps}=\rho_{\sf} \omega^2 \uu_{\eps}\cdot \nn+\sg 
%- \partial_{\nn }\sp_{i}
\quad&\mbox{on}\quad \Gamma
\\[0.5ex]
\TT(\uu_{\eps}) = -\sp_{\eps} \nn
%-\sp_{i} \nn 
 \quad &\mbox{on} \quad \Gamma
\\[0.2ex]
%\mbox{external B.C.} 
\sp_{\eps}=0
  \quad &\mbox{on}\quad \Gamma^\eps \ . % \partial_{\nn }\sp_{\eps} - i \kappa\sp_{\eps}=0 or Dirichlet
   \end{array}
    \right.
\end{equation}
%with a \Rd Dirichlet B.C. \Bk 
Hereafter, we explicit a weak formulation of the problem \eqref{EAfgh}. 
%with right-hand sides $f,g,h$ and
% with Dirichlet b.c. for  $(\uu,\sp)$ in the space
We first introduce the functional space adapted to a variational formulation
$$
V_{\eps}=
   \{ (\uu, \sp)  \in \bH^1(\Omega_{\mathsf{s}})\times \H^1 ({\Omega^{\eps}_\mathsf{f}}) \, |\, \gamma_{0} \sp =0\quad \mbox{on} \quad 
{\Gamma^\eps} \,
\} \ .
$$
Here, $\gamma_{0} \sp$ is the Dirichlet trace of $\sp $ on $\Gamma^\eps$. The space $V_{\eps}$ is endowed with the piecewise $\H^1$ norm in $\Omega_{\s}$ and $\Omega^{\eps}_\sf$. Then the  variational problem writes : Find $(\uu_{\eps},\sp_{\eps}) \in V_{\eps}$ such that
\begin{equation}
   \forall (\vv,\sq)  \in  V_{\eps}, \quad 
   a_{\eps}\left((\uu_{\eps},\sp_{\eps}), (\vv,\sq)\right) = 
   \left< F , (\vv,\sq) \right>_{V_{\eps}^{\prime},V_{\eps}}
 %  \overline{\Omega^\eps   }
  \ ,
\label{VP}
\end{equation}
where the sesquilinear form $a_{\eps}$ is defined as
\begin{multline*}
  a_{\eps} \left((\uu,\sp), (\vv,\sq)\right)=
   \int_{{\Omega^{\eps}_\mathsf{f}}} \left( \nabla \sp \cdot \nabla\bar\sq  -\kappa ^2 \sp \bar\sq \right)\,\dr\xx 
 +  \int_{\Omega_{\mathsf{s}}}\left(  \underline{\underline{\sigma}}( \uu) : \underline{\underline{\epsilon}} (\bar{\vv})  -\omega^2 \rho \uu  \bar{\vv}  \right)\dr\xx 
 \\
    + \int_{\Gamma}\left( \omega^2 \rho_{\sf} \uu\cdot \nn \, \bar\sq  + \sp  \bar{\vv} \cdot \nn \right)\,\dr \sigma \ ,
\end{multline*}
and the right-hand side $F$ is defined as
$$
%b\left( (\vv,\sq) \right)   
   \left< F , (\vv,\sq) \right>_{V_{\eps}^{\prime},V_{\eps}}
   %_{\overline{\Omega^\eps}} 
   =  -\int_{\Omega_\mathsf{f}^\eps} \sf   \bar\sq \, \dr\xx
    -\int_{\Omega_\s} \ff \cdot  \bar\vv\, \dr\xx
   - \int_{\Gamma} \sg \bar\sq  \,\dr \sigma \   .
 %  =   \int_{\Gamma}  \partial_{\nn} \sp_{i} \sq  - \sp_{i} \vv\cdot \nn   \,\dr s \   .
$$  
We assume that the data $(\ff,\sf,\sg)$ are smooth enough such that the right-hand side $F$ belongs to the space $V_{\eps}^{\prime}$.

\subsubsection*{Statement of uniform estimates}
In the framework of Sec. \ref{Sa} %(Assumptions \ref{AC}-\ref{SA}-\ref{HG}) 
we prove $\eps$-uniform a priori estimates for the solution of problem \eqref{VP}. The following theorem is the main result in this section.

\begin{thm}
\label{tUEeps}
Under Assumptions \ref{AC}-\ref{SA}-\ref{HG}, there exists constants $\eps_{0},C>0$ such that for all $\eps \in (0, \eps_{0})$, the problem \eqref{VP} with data  $F\in V_{\eps}^{\prime}$ has a unique solution  $(\uu_{\eps},\sp_{\eps}) \in V_{\eps}$ which satisfies % the uniform estimates
\begin{equation}
\label{UEeps}
%\sqrt{\eps}\| \partial_{t}  \fp_{\eps}\|_{0,\Omega_{\f}}+ \sqrt{\eps}^{-1} \| \partial_{3}  \fp_{\eps}\|_{0,\Omega_{\f}}+\| \fp_{\eps}\|_{0,\Omega_{\f}}
%+ \|  \fp_{\eps}\|_{0,\Gamma}
 \| \sp_{\eps}\|_{1,\Omega^\eps_{\sf}} 
+ \| \uu_{\eps}\|_{1,\Omega_{\s}} %+  \|\uu_{\eps}\cdot\nn\|_{0,\Gamma}
\leqslant C\| F\|_{V_{\eps}^{\prime}}\ .
 \end{equation}
\end{thm}
This result is proved in Sec. \ref{Sue}. The proof is based on a formulation of the problem set in a fixed domain. This formulation is obtained through a scaling along the thickness of the layer. %We adapt the proof of \cite{BL96}.
As an application of uniform estimates \eqref{UEeps}, we prove the convergence result \eqref{ek}  in Sec. \ref{SAEC}. %This result means that the Equivalent Condition in the asymptotic model \eqref{Eec} is of order $k$. %to solutions $(\uu_{\eps},\sp_{\eps})$ of \eqref{EA} as $\eps\to 0$.

\section{Equivalent Conditions}
\label{Smr}

 In the framework above, we derive for all $k\in\{0,1,2,3\}$  a boundary condition set on $\Gamma$  which is associated with the problem \eqref{EA} and satisfied by $\uu_\eps^{k}$, i.e. $\uu_\eps^{k}$ solves the problem
 \begin{equation}
\label{Ebck}
 \left\{
   \begin{array}{lll}
 \nabla \cdot  \underline{\underline{\sigma}} (\uu_\eps^{k})+\omega^2 \rho \uu_\eps^{k} =\ff
 \quad&\mbox{in}\quad \Omega_{\mathsf{s}}
\\[0.8ex]
\TT(\uu_\eps^{k})+\B_{k,\eps}(\uu_\eps^{k}\cdot\nn)\, \nn =0% \quad \mbox{and} \quad \uu \cdot\nn=0
  \quad &\mbox{on} \quad \Gamma\ .
 %\ 
   \end{array}
    \right.
\end{equation}
Here $\B_{k,\eps}$ is a surfacic differential operator acting on functions defined on $\Gamma$ and which depends on $\eps$. In this section, we present Equivalent Conditions (ECs) up to the fourth order and asymptotic models for the solution of the exact problem, Sec. \ref{S3DECs}. Then, we present well-posedness and convergence results, Sec. \ref{ScEc}.  Elements of derivation and mathematical validations for ECs are presented in Sec. \ref{SDEC} and Sec. \ref{SAEC}.

\subsection{Statement of  Equivalent conditions}
\label{S3DECs}
%\subsubsection{Dirichlet b.c}
We obtain a hierarchy of boundary-value problems. Each one gives a model with a different order of accuracy in $\eps$ and reflects the effect of the thin layer on the elastic displacement. We derive in Section \ref{CECs} the following boundary conditions in problem \eqref{Ebck} : 
\subsubsection*{Order 1} 
$$
\TT(\uu_{0}) =  0
%-\sp_{i} \nn  
\quad \mbox{on} \quad \Gamma\ 
$$
\subsubsection*{Order 2}
\begin{equation}
\label{Ebc1}
\TT(\uu_\eps^{1}) -\eps   \omega^2 \rho_{\sf} \uu_\eps^{1}\cdot \nn \, \nn= 0  
%-\sp_{i} \nn   -\eps    \partial_{\nn }\sp_{i} \nn 
 \quad \mbox{on} \quad \Gamma \ 
 \end{equation}
\subsubsection*{Order 3}
\begin{equation}
\label{Ebc2}
\TT(\uu_\eps^{2})- \eps   \omega^2 \rho_{\sf} \left( 1 -\eps \cH(y_\alpha) \right) \uu_\eps^{2}\cdot \nn \, \nn = 0 \quad \mbox{on} \quad \Gamma \ 
\end{equation}
%with $\hh_{2,\eps}= -\sp_{i} \nn -\eps   \left(1 + \eps  \cH(y_\alpha) \right)  \partial_{\nn }\sp_{i} \nn$.
\subsubsection*{Order 4}
\begin{multline}
\label{Ebc3}
\TT(\uu_\eps^{3}) -\eps   \omega^2 \rho_{\sf}\left(1 -\eps \cH(y_\alpha) +\frac{\eps^2}{3} \left[ \Delta_{\Gamma} + \kappa^2 \Id  + 4 \cH^2 (y_\alpha) -\cK(y_\alpha) \right]  \right)(\uu_\eps^{3}\cdot \nn) \, \nn
\\
=0  \quad \mbox{on} \quad \Gamma \ 
\end{multline}

Here, $(y_\alpha)$, $\alpha=1,2$, are tangential coordinates on $\Gamma$,  
%\[
%  \cH=\tfrac12\, b_{\alpha}^{\alpha}
%\]
$\cH$ and  $\cK$ denote the \textit{mean curvature} and the \textit{Gaussian  curvature} of the surface $\Gamma$ %, see Sec. \ref{SFAE},
 and $\Delta_{\Gamma}$ is the Laplace-Beltrami operator along $\Gamma$. Successive corrections appear in these conditions when increasing the order of approximation. The conditions of order  $k\in\{1,2,3\}$ involves only partial derivatives of order $2$ in the operator $\TT$, whereas the condition of order $k=4$ is a Ventcel's condition~\cite{AS78,Le87} since it involves partial derivatives of order $2$.

\begin{rem}
%In the equivalent condition of order $1$,
The ``background'' solution $\uu_{0}$ is independent of $\eps$.  It corresponds to a model where the thin layer is neglected. The effect of the fluid part appears from the order $2$ through the density $\rho_{\sf}$. The influence of the geometry of the surface $\Gamma$ appears from the order $3$ through the \textit{mean curvature} of $\Gamma$. The Helmholtz operator set on $\Gamma$ appears with the fourth order.  
\end{rem}

\begin{rem}
\iti1 (Comparison with~\cite{BL96}).   
We have compared impedance operators for elasto-acoustics and acoustics~\cite{BL96} precisely, Sec. \ref{CEca}. The comparison does not seem relevant.

\iti2 (Non-constant thickness). We consider here a thin layer with constant thickness.
In the context of geophysical applications, the thickness of the layer is no longer constant with respect to the tangential variable. The change of variables (or scaling) would lead to additional terms in the transmission
conditions. These terms come from the determinant of the metric of the layer. The derivation of the asymptotics are more tedious but all the tools are given in the present paper to perform the calculation.
\end{rem}

\subsection{Stability and convergence of Equivalent conditions}
\label{ScEc}

Our goal in the next sections is to validate ECs set on $\Gamma$ (Sec. \ref{S3DECs}) proving estimates for $\uu_{\eps}-\uu^k_{\eps}$ for all $k\in\{0,1,2,3\}$, where $\uu^k_{\eps}$ is the solution of the approximate model \eqref{Ebck}, and $\uu_{\eps}$ satisfies the problem \eqref{EA}. The functional setting for $\uu^k_{\eps}$ is described by the Hilbert space  $\bV^{k}$ :
\begin{notn}
\label{not2}
$\bV^{k}$ denotes the space $\bH^1(\Omega_{\mathsf{s}})$ when $k\in\{0,1,2\}$, and  $\{\uu\in \bH^1(\Omega_{\mathsf{s}}) \, | \quad \uu\cdot \nn|_{\Gamma} \in\H^1(\Gamma) \}$ when $k=3$\ .
\end{notn}

\begin{thm}
\label{CVec}
Under Assumptions \ref{AC}-\ref{SA}-\ref{HG}-\ref{Hspi}, for all $k\in\{0,1,2,3\}$ there exists constants $\eps_{k}, C_{k}>0$ such that for all $ \eps \in (0, \eps_{k})$, the  problem \eqref{Ebck} with data  $\ff\in \bL^2(\Omega_{\s})$ has a unique solution  $\uu^k_{\eps}\in \bV^{k}$ which satisfies uniform estimates
\begin{equation}
\label{UEec}
%\sqrt{\eps}\| \partial_{t}  \fp_{\eps}\|_{0,\Omega_{\f}}+ \sqrt{\eps}^{-1} \| \partial_{3}  \fp_{\eps}\|_{0,\Omega_{\f}}+\| \fp_{\eps}\|_{0,\Omega_{\f}}
%+ \|  \fp_{\eps}\|_{0,\Gamma}
% \| \sp_{\eps}\|_{1,\Omega^\eps_{\sf}} 
%+
 \| \uu_{\eps}-\uu^k_{\eps}\|_{1,\Omega_{\s}} %+  \|\uu_{\eps}\cdot\nn\|_{0,\Gamma}
\leqslant C_{k}\eps^{k+1} \ .
%\| F\|_{V_{\eps}^{\prime}}\ .
 \end{equation}
\end{thm}

The well-posedness result for the problem \eqref{Ebck} is stated in Thm. \ref{P1} and is proved in Sec. \ref{SWp}. It appears nontrivial to work straightforwardly with the difference $\uu_{\eps}-\uu^k_{\eps}$. 
%for a similar context in electromagnetism. 
A usual method consists to use the truncated series $\uu_{k,\eps}$ introduced in Sec. \ref{CECs} as intermediate quantities~\cite {HJN08}. Then, the error analysis is splitted into two steps detailed in the next sections : 
\medskip
\begin{enumerate}
\item We prove uniform estimates for the difference $\uu_{\eps}-\uu_{k,\eps}$ in Thm. \ref{TErem},  Sec. \ref{SEER} , 
\item We prove uniform estimates for the difference $\uu_{k,\eps}-\uu^k_{\eps}$, Sec. \ref{SEe} .  %Sec. \ref{Sii}. 
\end{enumerate}

%\begin{rem}
%The first step of the proof is independent of ECs and is valid for any integer $k$. The second step for $k=0$ is useless since $\uu_{0,\eps}=\uu^0_{\eps}$.
%\end{rem}

\subsection{Comparison with Equivalent conditions for acoustics}
\label{CEca}

In this section, we compare Equivalent conditions given in Sec. \ref{S3DECs} with the results for acoustics proved by Bendali--Lemrabet ~\cite{BL96}. ECs  for elasto-acoustics %in Sec. \ref{S3DECs} 
%are distinct from those given in  \cite{BL96}, since  they
 are of ``$\uu\cdot\nn$--to--$\TT(\uu)$'' nature, whereas ECs for acoustics involve local impedance operators which link the Dirichlet and Neumann traces. However, we eventually ``identify''  impedance operators involved  in Sec. \ref{S3DECs} and in~\cite{BL96}  through a Taylor expansion of the operator. 

\begin{figure}[ht]
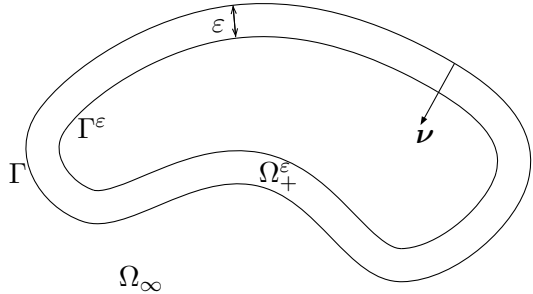

%%%%%%%%%%%%%%%%%%%%%%%%%%%%%%%%
%%%%%%%% DEBUT DESSIN
%%%%%%%%%%%%%%%%%%%%%%%%%%%%%%%
\def\EpsTube{-11}
\begin{center}
\psset(width=2)
\figinit{1.14pt}
\pssetupdate{yes}
% pt de la courbe \Sigma?
\figpt 8:(-48,-18)%\figpt 9:(-20,50)\figpt 10:(50,20)
%\figpt 12:(10,10)\figpt 13:(0,-10)
%\figpt 14:(-20,15)%
% pt de la courbe \Sigma_{h}
\figpt 22:(-92,25.5)\figpt 23:(-27.2,60.8)
%\figpt 23:(-37.2,30.8)
\figpt 24:(46.3,41.5)\figpt 25:(70.8,2.8)
\figpt 26:(29,-31)
%\figpt 27:(-27.2,8)
\figpt 27:(-20,1)\figpt 28:(-77,-11)
% pt de la courbe \Sigma
\figpt 1:(-26.3,50.2)
\figpt 2:(-12,-4.5)
\figpt 3:(50,50)
%\nn
\figpt 4:(35,21)\figpt 5:(30,-40)
%(46.3,41.5)
%\figpt 6:(-30,0)\figpt 7:(-80,-20)
% ecriture de \Sigma, \Omega\iso, \Omega\con et \Partial\Omega, et \nn et h et ${\mathcal{O}}$
\figpt 16:(-70,20)\figpt 15:(0,30)
\figpt 17:(40,10) \figpt 18:(-60,-20)
\figpt 21:(50,65) \figpt 31:(30,-18)  \figpt 34:(-32,65) 
 \figpt 35:(-89,20)
   \figpt 36:(-103,5)
   \figpt 37:(-70,-30)
% pt vecteurs normaux
\figpt 19:(61,65) \figpt 20:(3,27)  \figpt 30:(42.6,33)  \figpt 32:(55,7) \figpt 33:(61.6,6.1) 
%%%%%%%%%%%%%%%%%%%%%
%% fichier graphique
%%%%%%%%%%%%%%%%%%%
\psbeginfig{}
\psset arrowhead(fillmode=no,length=3)\psarrow[1,23]
\psset arrowhead(fillmode=no,length=3)\psarrow[23,1]
%\psarrow[23,2]
%\pssetfillmode{yes}\pssetgray{0.5}
%\pscurve[1,2,3,4,5,6,7,1,2,3]
%\psBezier 4[1,2,3,4,5,6,7,1]
%\pscurve[8,9,10,12,13,14,8,9,10]
\pssetfillmode{no}\pssetgray{0}
%\pscurve[8,9,10,12,13,14,8,9,10]
%\pssetfillmode{yes}\pssetgray{0.5}
\pscurve[22,23,24,25,26,27,28,22,23,24]
\figptscontrolcurve 40,\NbC[22,23,24,25,26,27,28,22,23,24]
%\psreset{first}\psset (width=1.9)
\psEpsLayer \EpsTube,\NbC[40,41,42,43,44,45,46,47,48,49,50,51,52,53,54,55,56,57,58,59,60,61]
%\pscurve[40,43,46,49,52,55,58,40,43]
%\figptscontrolcurve 70,\NbC[40,43,46,49,52,55,58,40,43]
%\psreset{first}\psset (width=1.9)
%\psEpsLayer \EpsTubebis,\NbC[70,71,72,73,74,75,76,77,78,79,80,81,82,83,84,85,86,87,88,89,90,91]
%\pscurve[22,23,24,25,26,27,28,22,23,24]
%\figptscontrolcurve 70,\NbC[22,23,24,25,26,27,28,22,23,24]
%\psreset{first}\psset (width=1.9)
%\psEpsLayer \EpsTubebis,\NbC[70,71,72,73,74,75,76,77,78,79,80,81,82,83,84,85,86,87,88,89,90,91]
%\'epaisseur h
%\psline[23,2]
%\psarrow[23,2]
%\psarrow[2,23]
%\pssetfillmode{yes}\pssetgray{0.5}
%Vecteurs
%\psarrow[32,25]
\psset arrowhead(fillmode=yes,length=3)\psarrow[24,4]
\pssetfillmode{yes}\pssetgray{0.5}
\psendfig
%%%%%%%%%%%%%%%%%%%
%% writing
%%%%%%%%%%%%%%%%%%%
\figvisu{\figBoxA}{
%{Figure 1}\ --\ Un voisinage tubulaire $\mathcal{O^{\prime}}\subset\mathcal{O}$ de $\Sigma$ dans $\Omega$
}
{
\figwritee 37: $\Omega_\infty$ (2pt)
%\figwritec [16]{$\Omega\iso$}
\figwriten 2: $\Omega^{\eps}_+$ (4pt)
\figwrites 34 : $\eps$ (10pt)
\figwritee 36: $\Gamma$ (-2pt)
%\figwriten 31: $\Sigma_{h}$ (2pt)
\figwrites 4: $\nnu$(2pt)
%\figwritew 32: $\nn$(3pt)
%\figwritee 20: $\nn$(3pt)
\figwritee 35 : $\Gamma^\eps$ (8pt)
\figsetmark{$\figBullet$}
%\figwritep[1,2,3,4,5,6,7,22,23,24,25,26,27,28]
}
\centerline{\box\figBoxA}
 \caption{
 %A cross-section of the domain $\Omega^\eps$ and its subdomains $\Omega_\mathsf{s}$ and $\Omega^{\eps}_\mathsf{f}$ 
 Geometry of the considered problem \eqref{EABL}
 }
\label{F1}
\end{center}
\end{figure}

 We first remind the results of the works in Ref.~\cite{BL96} where the authors consider the scalar problem in time-harmonic regime
\begin{equation}
\label{EABL}
 \left\{
   \begin{array}{lll}
  \Delta \su^+_{\eps}+\kappa^2  n^2 \su^+_{\eps} = 0   \quad&\mbox{in}\quad \Omega^{\eps}_+
\\[0.5ex]
  \Delta \su^-_{\eps}+\kappa^2  \su^-_{\eps} = \sf   \quad&\mbox{in}\quad \Omega_{\infty}
\\[0.5ex]
\partial_{\nnu}\su^+_{\eps}= \alpha^{-1} \partial_{\nnu}\su^-_{\eps} \quad&\mbox{on}\quad \Gamma
\\[0.5ex]
\su^+_{\eps}= \su^-_{\eps} \quad &\mbox{on} \quad \Gamma
\\[0.2ex]
\su^+_{\eps}=0   \quad &\mbox{on}\quad \Gamma^\eps     
\\[0.2ex]
\mbox{Sommerfeld b.c.} \quad &\mbox{at}\quad  \infty \ ,    
\end{array}
    \right.
\end{equation}
set in a smooth unbounded domain % $\Omega^\eps$
 in $\R^3$ made of    a smooth connected subdomain $\Omega_{\infty}$ and a subdomain $\Omega^{\eps}_+$ which is a thin layer of uniform thickness $\eps$,  see figure \ref{F1}. Here, $n>0$ represents a refractive index, $\alpha>0$ and  $\sf$ is a data. 
In this framework, the authors derive in~\cite{BL96} equivalent conditions set on $\Gamma$ which are associated with the problem \eqref{EABL} and satisfied by $\su_\eps^{k}$ for all $k\in\{0,1,2,3\}$, i.e. $\su_\eps^{k}$ solves the problem
\begin{equation}
\label{EbckBL}
 \left\{
   \begin{array}{lll}
  \Delta \su_\eps^{k}+\kappa^2 n^2 \su_\eps^{k} = \sf   \quad&\mbox{in}\quad \Omega_{\infty}\ 
\\[0.5ex]
\sN_{\eps,k}(\su_\eps^{k})+\partial_{\nnu}\su_\eps^{k} =0  \quad &\mbox{on} \quad \Gamma\ 
\\[0.2ex]
\mbox{Sommerfeld b.c.} \quad &\mbox{at}\quad  \infty \ ,
   \end{array}
    \right.
\end{equation}
Equivalent conditions in \eqref{EbckBL} write successively for $k\in\{0,1,2,3\}$:
%The impedance operators $\sN_{\eps,k}$ in {EbckBL} are local :

\subsubsection*{Order 1} 
$$
\su_{0}=  0 \quad \mbox{on} \quad \Gamma\ 
$$
\subsubsection*{Order 2}
$$
 \alpha \eps^{-1}  \su_\eps^{1} +\partial_{\nnu} \su_\eps^{1}= 0  
 \quad \mbox{on} \quad \Gamma \ 
$$
\subsubsection*{Order 3}
$$
 \alpha \eps^{-1} \left( 1 + \eps \cH(y_\alpha) \right)  \su_\eps^{2} + \partial_{\nnu} \su_\eps^{2}= 0 \quad \mbox{on} \quad \Gamma \ 
$$
\subsubsection*{Order 4}
$$
 \alpha \eps^{-1}
  \left( 1 + \eps \cH(y_\alpha) +\frac{\eps^2}{3} 
  \left[- \Delta_{\Gamma} - \kappa^2 n^2 \Id  +\cK(y_\alpha)-  \cH^2 (y_\alpha)  \right] \right)  \su_\eps^{3}  
 +  \partial_{\nnu} \su_\eps^{3}=0  \quad \mbox{on} \quad \Gamma \ 
$$

Setting $\alpha^{-1}=\rho \omega_{\sf}$, $n=1$, and using a formal Taylor expansion of the operator $\sN^{-1}_{\eps,k}$ ($k\in\{2,3\}$), there holds
\begin{subequations}
\begin{gather}
\label{sNB1}
\sN^{-1}_{\eps,1}= -\B_{\eps,1}
 \, \\
 \label{sNBk}
\sN^{-1}_{\eps,k}= -\B_{\eps,k}+\mathcal{O}(\eps^{k+1})  \quad \mbox{when} \quad k=2,3 
\, .
\end{gather}
\end{subequations}
%Hence, we can "identify" 
 %through 
 The relations \eqref{sNB1}-\eqref{sNBk} make the links between 
 the operators $\B_{\eps,k}$ and $\sN_{\eps,k}$ involved in the ECs  \eqref{Ebck} and \eqref{EbckBL}, respectively.
 % the ``$\uu\cdot\nn$--to--$\TT(\uu)$'' conditions \eqref{Ebck} and in the Dirichlet--to--Neumann conditions \eqref{EbckBL}. %through the relations \eqref{sNB1}-\eqref{sNBk}. 

\section{Uniform Estimates}
\label{Sue}

In this section, we prove uniform estimates for the exact solution of the elasto-acoustic problem. Since the functional setting of the variational problem \eqref{VP} depends on the small parameter $\eps$, it is not well suited to prove uniform estimates for solutions $(\uu_{\eps},\sp_{\eps}) \in V_{\eps}$. To overcome this difficulty, we adapt an idea developed  in~\cite {BL96} writing equivalently the variational problem \eqref{VP} in a common functional framework as $\eps$ varying, Sec. \ref{Ssp}.  We state uniform estimates in this new framework, Th. \ref{UE}. We use a compactness argument to prove estimates, Sec. \ref{S4lem}.

\subsection{The scaled problem}
\label{Ssp}

%Following \cite[\S 4.2]{Pe13}, 
We write the variational problem \eqref{VP} in a fixed domain through the scaling $S=\eps^{-1} \nu $ where $\nu \in(0,\eps)$ is the distance to the surface $\Gamma$. The fixed domain writes $\Omega_{\mathsf{s}}\times\Omega_\f$ where $\Omega_\f:= \Gamma\times (0,1)$ and the ad-hoc functional space writes 
$$
V= \{ (\uu, \fp)  \in \bH^1(\Omega_{\mathsf{s}})\times \H^1 (\Omega_\f) \, |\quad \fp(.\,,1)=0\quad \mbox{on} \quad \Gamma 
\} \ .
$$
Then the variational problem writes :
Find $(\uu_{\eps},\fp_{\eps}) \in V$ such that for all $(\vv,\fq)  \in  V$ , 
\begin{equation}
    \eps  a_{\sf}(\eps; \fp_{\eps},\fq)+
        a_{\s}( \uu_{\eps},\vv)
         +
        \int_{\Gamma}\left( \omega^2 \rho_{\sf} \uu\cdot \nn \, \bar\fq  + \fp  \bar{\vv} \cdot \nn \right)\,\dr \Gamma
        = 
     \left< \F_{\eps} , (\vv,\fq) \right>_{V' , V}  \ ,
\label{VFinV}
\end{equation}
where 
%\begin{multline}
\begin{equation*}
  a_{\sf}(\eps; \fp,\fq)=
 \int_{0}^1\int_{\Gamma} 
  \left\{ 
 \left(\Id+\eps S \mathcal{R}\right)^{-2}  \nabla_{\Gamma}\fp  \nabla_{\Gamma}\bar\fq
 +\eps^{-2} \partial_{S} \fp  \partial_{S} \bar \fq 
   -\kappa^2  \fp \bar \fq
 \right\}  \det \left(\Id+\eps S \mathcal{R}\right) \dr \Gamma \dr S
%  \\
 % -\kappa^2 \int_{0}^1\int_{\Gamma} \fp \bar \fq \, \det \left(\Id+\eps S \mathcal{R}\right)   \dr \Gamma \dr S
  \ ,
  \end{equation*}
%  In \eqref{VFinV}, 
\begin{equation*}
     a_{\s}( \uu,\vv)=   \int_{\Omega_{\mathsf{s}}}\left( \underline{\underline{\sigma}} (\uu) : \underline{\underline{\epsilon}} (\bar{\vv})  -\omega^2 \rho \uu  \bar{\vv}  \right)\dr\xx \ ,
     %\\
%+ \int_{\Omega^{\eps}_\mathsf{f}} \left( \nabla \sp \cdot \nabla\bar\sq  -\kappa ^2 \sp \bar\sq \right)\,\dr\xx \\
  %  + \int_{\Gamma}\left( \omega^2 \rho_{\sf} \uu\cdot \nn \, \bar\sq  + \sp  \bar{\vv} \cdot \nn \right)\,\dr \Gamma \ ,
\end{equation*}
\begin{equation*}
\quad \mbox{and} \quad   \left< \F_{\eps} , (\vv,\fq) \right>_{V' , V}  =  -\eps\int_{\Omega_\f} \f   \bar\fq  
   \det \left(\Id+\eps S\mathcal{R}\right) \dr \Gamma \dr S 
   -\int_{\Omega_\s} \ff \cdot  \bar\vv\, \dr\xx
     - \int_{\Gamma}  \fg \bar\fq  \,\dr \Gamma \ .
\end{equation*}
%\end{multline}
Here $\mathcal{R}$ is an intrinsic symmetric linear operator  defined on the tangent plane $\TT_{\xx_{\Gamma}}(\Gamma)$ to $\Gamma$ at the point $ \xx_{\Gamma}\in\Gamma$ which characterizes the curvature of $\Gamma$ at the point $ \xx_{\Gamma}$. We refer to~\cite{BL96}-\cite[\S 4.1]{Pe13} for the introduction of geometrical tools and more details. The parameter  $\eps$ weighting the form $a_{\sf}(\eps; \fp,\fq)$ in formulation \eqref{VFinV} may lead to a solution  $(\uu_{\eps},\fp_{\eps}) \in V$ such that the surface gradient $\nabla_{\Gamma} \fp_{\eps}$ can be unbounded as $\eps\to 0$. This is a similarity with the work in Ref.~\cite {BL96}. 
 %We refer to~\cite {BL96} where the authors exhibit this kind of singular perturbation term for a diffraction problem of acoustic waves.
 Furthermore, the sign of the left-hand side of the problem \eqref{VFinV}  for $(\vv, \fq)=(\uu_{\eps},\fp_{\eps})$ cannot be controlled. 
 %As $\nabla_{\Gamma} \fp_{\eps} $ may be unbounded, it can generate a solution $(\uu_{\eps},\fp_{\eps})$ depending singularly of $\eps$ too in $\Omega_{\s}$ for general right-hand sides. Fortunately, the right-hand sides which appear in the asymptotic analysis of problem  \eqref{VFinV} compensate the singular behavior of the perturbation term.  
Hence due to the lack of strong coerciveness of the variational formulation  \eqref{VFinV}  one cannot 
%take advantage of the previous fact to
 get  straightforwardly estimates. %: the formulation is not sufficient by itself to ensure the stability with respect to $\eps$ on which the asymptotic analysis is based. Therefore the proof of this stability result constitutes a main part of the asymptotic analysis.
 % The proof of this result involves both a compactness  argument and the spectral assumption \ref{SA}. 
%of the related problem obtained when completely neglecting the effect of the layer. %Thus, the amount of smoothness of $\Gamma$ is an essential ingredient in the asymptotic analysis 
%We first prove some preliminary results (Poincar\'e and trace inequalities) useful for the proof of Lemma \ref{L2}.
 % Using the projection operator $\Pi_{ \xx_{\Gamma}}$ from $\R^3$ onto the tangent plane $\TT_{ \xx_{\Gamma}}(\Gamma)$, there holds : 
%\begin{equation*}
%\label{Enablasp}
%\Pi_{ \xx_{\Gamma}} \nabla \sp=  \left(\Id +\eps S \mathcal{R} \right)^{-1}\nabla_{\Gamma} \fp \ , %\quad [36] 
 %\end{equation*}
Our main result for the problem \eqref{VFinV} is the following a priori estimate, uniform as $\eps\to 0$.

\begin{thm}
\label{UE}
Under Assumptions \ref{AC}-\ref{SA}-\ref{HG},  there exists constants $\eps_{0},C>0$ such that for all $ \eps \in (0, \eps_{0})$, the problem \eqref{VFinV} with data  $\F_{\eps}\in V^{\prime}$ has a unique solution  $(\uu_{\eps},\fp_{\eps}) \in V$ which satisfies  the uniform estimates
\begin{equation}
\label{UEup}
\sqrt{\eps}\| \nabla_{\Gamma}  \fp_{\eps}\|_{0,\Omega_{\f}}+ \sqrt{\eps}^{-1} \| \partial_{S}  \fp_{\eps}\|_{0,\Omega_{\f}}+\| \fp_{\eps}\|_{0,\Omega_{\f}}
+ \|  \fp_{\eps}\|_{0,\Gamma}
+ \| \uu_{\eps}\|_{1,\Omega_{\s}} %+  \|\uu_{\eps}\cdot\nn\|_{0,\Gamma}
\leqslant C\| \F_{\eps}\|_{V^{\prime}}\ .
 \end{equation}
\end{thm}
This theorem is the key for the proof of Thm. \ref{tUEeps} : as a consequence of the following estimates  %\eqref{Ea}-\eqref{Eb} 
\begin{subequations}
\begin{gather}
\label{Ea}
\forall \sp\in \L^2(\Omega^\eps_{\sf})\ , \quad  \| \sp \|_{0,\Omega^\eps_{\sf}}\simeq \sqrt{\eps}   \| \fp \|_{0,\Omega_{\f}}
 \,,\\
\label{Eb}
\forall \sp\in \H^1(\Omega^\eps_{\sf})\ , \quad    \|\nabla \sp \|_{0,\Omega^\eps_{\sf}}\simeq \sqrt{\eps}   \| \nabla_{\Gamma}\fp \|_{0,\Omega_{\f}}+ \sqrt{\eps}^{-1}   \| \partial_{S}\fp \|_{0,\Omega_{\f}}
\, ,
\end{gather}
\end{subequations}
available  for $\eps$ small enough, and a proof of which can be found in~\cite[\S 4.1]{Pe13}, we obtain estimates  \eqref{UEeps}. 
% and we infer the Thm. \ref{tUEeps}.
In \eqref{Ea}-\eqref{Eb},  for any function $\sp$ defined in $ \Omega^{\eps}_\mathsf{f}$, the function $\fp$ is defined in the domain  $\Omega_\f$ as 
$$
  \fp( \xx_{\Gamma},S)  = \sp(\xx)\ , \quad ( \xx_{\Gamma},S=\frac{\nu}{\eps})\in \Gamma \times (0,1) \ .
$$
In \eqref{Ea}, the symbol $\simeq$ means that quantities $\| \sp \|_{0,\Omega^\eps_{\sf}}$ and $\sqrt{\eps}   \| \fp \|_{0,\Omega_{\f}}$ are equal up to a multiplicative constant which is independent of $\eps$.

The proof of Thm. \ref{UE} is based on the following statement.
%We consider now the scaled problem \eqref{VFinV}  at a fixed frequency $\omega$ satisfying Assumption \ref{SA}. We are going to prove the following statement : 

\begin{lem}
\label{L2}
Under Assumptions \ref{AC}-\ref{SA}-\ref{HG}, there exists  constants $\eps_{0},C>0$ such that for all $ \eps \in (0, \eps_{0})$, any solution $(\uu_{\eps},\fp_{\eps}) \in V$ of problem \eqref{VFinV} with a data $\F_{\eps}\in V^{\prime}$ satisfies the uniform estimates
\begin{equation}
\label{UEL2up}
 \| (\uu_{\eps},\fp_{\eps})\|_{0,\Omega_{\s}\times\Omega_{\f}} +  \|\uu_{\eps}\cdot\nn\|_{0,\Gamma}+
 \|  \fp_{\eps}\|_{0,\Gamma}\leqslant C\| \F_{\eps}\|_{V^{\prime}}\ .
 \end{equation}
\end{lem}

This Lemma is going to be proved in Sec. \ref{S4lem}.  The proof of this result involves both a compactness  argument and the spectral assumption \ref{SA}. 
 As a consequence of estimates \eqref{UEL2up}, we infer estimates  \eqref{UEup}. Since the problem \eqref{VFinV} is of Fredholm type, the Thm. \ref{UE} is then obtained as a consequence of the Fredholm alternative. %This argument is rather standard, we refer for instance to the proof of ~\cite[Cor.\,4.3]{CDP10} for a similar argument.

%\begin{lem}
%\label{L3}
%Under Assumption \ref{SA}, $\exists \eps_{0}>0$, $\forall \eps \in (0, \eps_{0})$,  for all $(\uu_{\eps},\fp_{\eps}) \in V\subset W$ solution of problem \eqref{VFinV} associated with $\F\in V^{\prime}$, there holds the uniform estimates
%\begin{equation}
%\label{UEup}
%\sqrt{\eps}\| \partial_{t}  \fp\|_{0,\Omega_{\f}}+ \sqrt{\eps}^{-1} \| \partial_{3}  \fp\|_{0,\Omega_{\f}}+\| \fp_{\eps}\|_{0,\Omega_{\f}}
%+ \|  \fp\|_{0,\Gamma}+
 %\| \uu_{\eps}\|_{1,\Omega_{\s}} %+  \|\uu_{\eps}\cdot\nn\|_{0,\Gamma}
%\leqslant C\| \F\|_{V^{\prime}}\ .
 %\end{equation}
%\end{lem}

\subsection{Proof of Lemma \ref{L2}\,: Uniform estimate of  $(\uu_{\eps},\fp_{\eps})$}
\label{S4lem}
We prove this lemma by contradiction : We assume that there exists a sequence $(\uu_m,\fp_{m})\in V$, $m\in\N$, of solutions of problem \eqref{VFinV} associated with a parameter $\eps_m$ and a right-hand side  $\F_{m}\in V^{\prime}$:
\begin{equation}
   \forall (\vv,\fq)  \in  V, \quad 
   \eps_{m}  a_{\sf}(\eps_{m}; \fp_{m},\fq)+
  %   \eps  b_{\sf}(\eps_{m}; \fp_{m},\fq)+
        a_{\s}( \uu_{m},\vv)
        +
        \int_{\Gamma}\left( \omega^2 \rho_{\sf} \uu_{m}\cdot \nn \, \bar\fq  + \fp_{m}  \bar{\vv} \cdot \nn \right)\,\dr \Gamma
        = 
   %b\left( (\vv,\sq) \right)
     \left< \F_{m} , (\vv,\fq) \right>_{V^{\prime},V} \ ,  
\label{VFinVm}
\end{equation}
%\begin{subequations}
%\begin{gather}
%\label{E3a}
% \omega^2 \rho \uu_m + \nabla \cdot  \underline{\underline{C}} \ \underline{\underline{\nabla}} \uu_m=0
 %\quad\mbox{in}\quad\Omega\,,\\
%\label{E3b}
 % \TT(\uu_{m}) -\eps_{m} \omega^2 \rho_{\sf} \uu_{m}\cdot \nn \, \nn=  \h_{m}  \nn\quad\mbox{on}\quad\Gamma\,,
%\end{gather}
%\end{subequations}
satisfying the following conditions
\begin{subequations}
\begin{eqnarray}
\label{5Ea}
   &\eps_m\to 0\quad &\mbox{as \ $m\to\infty$,} \\
\label{5Eb}
   & \| (\uu_{m},\fp_{m})\|_{0,\Omega_{\s}\times\Omega_{\f}} +  \|\uu_{m}\cdot\nn\|_{0,\Gamma}+
 \|  \fp_{m}\|_{0,\Gamma}
   %\|\uu_m\|_{0,\Omega} 
   = 1\quad &\mbox{for all $m\in\N$,} \\
\label{5Ec}
   &\|\F_m\|_{V^{\prime}}\to0\quad&\mbox{as \ $m\to\infty$.}
\end{eqnarray}
\end{subequations}

%We particularize the elastic variational formulation \eqref{FV} for the sequence $\{\uu_m\}$: For all $\vv\in\bH^1(\Omega)$:
%\begin{equation}
%\label{EVm}
 % \int_{\Omega_{\mathsf{s}}}\left(  \underline{\underline{C}} \ \underline{\underline{\epsilon}}  (\uu_{m}) : \underline{\underline{\epsilon}}  (\bar{\vv})  -\omega^2 \rho \uu_{m}  \bar{\vv}  \right)\;\dr\xx  \ -\eps_{m}\omega^2\rho_{\sf} \int_{\Gamma}    \uu_{m} \cdot \nn\ \bar{\vv} \cdot \nn \,\dr s 
%= \int_{\Gamma}  \h_{m}  \ \bar{\vv} \cdot\nn \,\dr s \   .
%\end{equation}
Choosing tests functions $(\vv,\fq)=(\uu_m,\fp_{m})$ in \eqref{VFinVm}, we obtain with the help of conditions \eqref{5Ea}-\eqref{5Eb}-\eqref{5Ec} the following uniform bounds : 

\iti 1 The sequence $\{\uu_m,\fp_{m}\}$ is bounded in the space $W$ defined as
 $$
W= \{ (\uu, \fp)  \in 
\bH^1(\Omega_{\mathsf{s}}) \times \H^1 (0,1; \L^2(\Gamma))
\, |\quad \fp(.\,,1)=0 \quad \mbox{on} \quad \Gamma \,\} \ ,
$$
% According to  \eqref{5Ec} we infer that  the sequence $\{\uu_m\}$ is bounded in $\bH^1(\Omega)$ :
\begin{equation}
\label{3E12}
   \|(\uu_m,\fp_{m})\|_{W}  
   \leqslant C .
\end{equation}
Remind that $\H^1 (0,1; \L^2(\Gamma))$ is the space of distributions   $\fp
\in\mathcal{D}^{\prime}(0,1; \L^2(\Gamma))$ such that $\fp$ and $\fp^{\prime}$ belong to $\L^2(0,1; \L^2(\Gamma))$. Subsequently, we identify the space $\L^2(0,1; \L^2(\Gamma))$ and $\L^2(\Omega_{\f})$. 

Since the domain $\Omega_{\s}$ is bounded, the embedding of $\bH^{1}(\Omega_{\s})$ into $\bL^2(\Omega_{\s})$ is compact.  
%and the embedding of $\bH^{\frac12}(\Gamma)$ into $\L^2(\Gamma)$ is compact
 Hence as a consequence of \eqref{3E12}, using the Rellich Lemma we can extract a subsequence  of $\{\uu_m,\fp_{m}\}$ (still denoted by $\{\uu_m,\fp_m\}$) which is converging in $\bL^2(\Omega_{\s})\times \L^2(\Omega_{\f})$, and we can assume that the sequence $\{ \underline{\underline{\nabla}} \uu_m\}$ is weakly converging in $\LL^2(\Omega_{\s})$.

\iti2 The sequence $\{\sqrt{\eps_{m}}^{-1} \partial_{S}\fp_{m}\} $ is bounded in $\L^2(\Omega_{\f})$, hence the sequence $\{ \partial_{S}\fp_{m}\}$ converges  to $0$ in 
$\L^2(\Omega_{\f})$.

\subsubsection*{Limit of the sequences $\{\uu_m,\fp_{m}\}$}
To prove convergence results for  $\{\uu_m,\fp_{m}\}$, we use Poincar\'e and trace inequalities which are available in the Hilbert space 
$$
\H= \{  \fp  \in  \H^1 (0,1; \L^2(\Gamma))
 \, |\quad \fp(.\,,1)=0 \quad \mbox{on} \quad \Gamma \,\} \ ,
 $$
see~\cite[Prop. 4.4 - 4.5]{Pe13} : There exists $C_{P},C>0$ such that 
\begin{equation}
\label{EP}
\forall \fp\in \H \ ,\quad 
\|  \fp\|_{0,\Omega_{\f}}\leqslant C_{P} \|  \partial_{S} \fp\|_{0,\Omega_{\f}}
 \quad \mbox{and} \quad
 \| \gamma_{0} \fp\|_{0,\Gamma}\leqslant C \|  \partial_{S} \fp\|_{0,\Omega_{\f}} \ .
\end{equation}
%For any $S\in (0,1)$, the trace $\fp(.\, ,S)$ is well defined since $S \mapsto  \fp(.\, ,S)$ is a continuous function of $S$ with values in $\L^2(\Gamma)$.

As a consequence of \iti2 and \eqref{EP}, we infer that  sequences $\{\fp_{m}\} $ and $\{\gamma_{0}\fp_{m}\}$ converge  to $0$ in 
$\L^2(\Omega_{\f})$ and 
%the sequence $\{\gamma_{0}\fp_{m}\}$  converges 
 to $0$ in $\L^2(\Gamma)$ respectively
\begin{equation}
\left\{
   \begin{array}{lll}
   \label{Ecvpm}
   \fp_m \rightarrow 0 \quad &\mbox{in} \quad \L^2(\Omega_{\f})
\\
   \gamma_{0} \fp_m \rightarrow 0   \quad &\mbox{in} \quad \L^2(\Gamma).
   \end{array}
\right.
\end{equation}
%There exists $C>0$ such that 
%and 
%$$
%\forall \fp\in \H \ , 
%\quad \|  \fp\|_{0,\Gamma}\leqslant C \|  \partial_{S} \fp\|_{0,\Omega_{\f}} 
%$$
%\end{prop}
%.

 Another consequence of  \eqref{3E12} is that the sequence $\{\uu_m \cdot\nn\}$ is bounded in $\H^{\frac12}(\Gamma)$. Therefore (up to the extraction of a subsequence) we can assume that the sequence $\{\uu_m \cdot\nn\}$ is strongly converging in $\L^2(\Gamma)$. Summarizing these convergence results, we deduce that there exists $\uu\in \bL^2(\Omega_{\s})$ such that
\begin{equation}
\left\{
   \begin{array}{lll}
   \label{3E6}
 \underline{\underline{\epsilon}} (\uu_m) \rightharpoonup   \underline{\underline{\epsilon}} (\uu)   \quad &\mbox{in} \quad \LL^2(\Omega_{\s})
\\
   \uu_m \rightarrow \uu   \quad &\mbox{in} \quad \bL^2(\Omega_{\s})
   \\
     \uu_m \cdot\nn \rightarrow \uu \cdot\nn  \quad &\mbox{in} \quad \L^2(\Gamma).
   \end{array}
\right.
\end{equation}

As a consequence of the strong convergence of sequences $\{\uu_{m}\}$ in $\bL^2(\Omega_{\s})$ and $\{\uu_m \cdot\nn\}$ in $\L^2(\Gamma)$, and  the strong convergence of $\fp_{m}$ and  $\gamma_{0}\fp_{m}$ \eqref{Ecvpm} together with \eqref{5Eb}, we infer 
$$
 \| \uu\|_{0,\Omega_{\s}} +  \|\uu\cdot\nn\|_{0,\Gamma}  = 1\ .
 $$

\subsubsection*{Conclusion}

Using Assumption \ref{SA}, we are going  to prove hereafter that $\uu=0$, which will contradict $\|\uu\|_{0,\Omega_{\s}}+  \|\uu\cdot\nn\|_{0,\Gamma}  =1$, and finally prove estimate \eqref{UEL2up}. We use $(\vv, \fq=0)\in V$ as test functions in \eqref{VFinVm}: there holds
 
\begin{equation*}
 %  \forall (\vv,\fq)  \in  V, \quad 
 %  \eps  a_{\sf}(\eps_{m}; \fp_{m},\fq)+
   %  \eps  b_{\sf}(\eps_{m}; \fp_{m},\fq)+
   %     a_{\s}( \uu_{m},\vv)
         \int_{\Omega_{\mathsf{s}}}\left( \underline{\underline{\sigma}}  (\uu_{m}) : \underline{\underline{\epsilon}}  (\bar{\vv})  -\omega^2 \rho \uu_{m}  \bar{\vv}  \right)\;\dr\xx
        \\
        +
        \int_{\Gamma}
        %\left( \omega^2 \rho_{\sf} \uu_{m}\cdot \nn \, \bar\fq  +
         \fp_{m}  \bar{\vv} \cdot \nn
         % \right)
         \,\dr \Gamma
        = 
   %b\left( (\vv,\sq) \right)
     \left< \F_{m} , (\vv,0) \right>_{V^{\prime},V} \ .
%\label{VFinVm}
\end{equation*}

%Taking the limit in both sides of the previous equation, we infer
%\begin{equation*}
%\label{EVm}
 % \int_{\Omega_{\mathsf{s}}}\left(  \underline{\underline{C}} \ \underline{\underline{\epsilon}}  (\uu_{m}) : \underline{\underline{\epsilon}}  (\bar{\vv})  -\omega^2 \rho \uu_{m}  \bar{\vv}  \right)\;\dr\xx  \ -\eps_{m}\omega^2\rho_{\sf} \int_{\Gamma}    \uu_{m} \cdot \nn\ \bar{\vv} \cdot \nn \,\dr s 
%= \int_{\Gamma}  \h_{m} \ \bar{\vv} \cdot\nn \,\dr s \   .
%\end{equation*}
According to \eqref{3E6}, \eqref{Ecvpm} and \eqref{5Ec}, taking limits as $m\rightarrow +\infty$, we deduce from the previous equalities $\uu \in \bH^1(\Omega_{\s})$ satisfies for all $\vv \in \bH^1(\Omega_{\s})$:
\begin{equation}
\label{E20}
  \int_{\Omega_{\mathsf{s}}}\left(   \underline{\underline{\sigma}}  (\uu) : \underline{\underline{\epsilon}}  (\bar{\vv})  -\omega^2 \rho \uu \bar{\vv}  \right)\;\dr\xx 
  % \ -\eps_{m}\omega^2\rho_{\sf} \int_{\Gamma}    \uu_{m} \cdot \nn\ \bar{\vv} \cdot \nn \,\dr s 
= 0\ .
%\int_{\Gamma}  \hh_{m} \ \bar{\vv} \cdot\nn \,\dr s \   .
\end{equation}
Integrating by parts the first term in the sesquilinear form \eqref{E20}, %there holds
%\begin{equation}
%\label{ipp}
%  - \left( \underline{\underline{\sigma}}( \uu) , \underline{\underline{\epsilon}}(\vv)\right)_{0,\Omega_{\mathsf{s}}} = 
%   \left( \nabla \cdot \underline{\underline{\sigma}}( \uu) , \vv\right)_{0,\Omega_{\mathsf{s}}} -\left(\TT(\uu),\vv\right)_{|\Gamma}.
%\end{equation}
we find that $\uu$ satisfies the problem
\begin{equation*}
 \left\{
   \begin{array}{lll}
\nabla \cdot \underline{\underline{\sigma}}( \uu) +\omega^2 \rho \uu =0
 \quad&\mbox{in}\quad \Omega_{\mathsf{s}}
\\[0.5ex]
\TT(\uu) =0 \quad &\mbox{on} \quad \Gamma \ .
   \end{array}
    \right.
\end{equation*}
 %Hence $\uu\in\bH^1(\Omega)$ is solution of problem \eqref{H}. 
 By Assumption \ref{SA}, we deduce
\begin{equation*}
  \uu=0 \quad \mbox{in} \quad \Omega_{\s}\ ,
\end{equation*}
which contradicts 
$\|\uu\|_{0,\Omega}+  \|\uu\cdot\nn\|_{0,\Gamma} =1$ and ends the proof of Lemma \ref{L2}.

\section{Derivation of Equivalent Conditions}
\label{SDEC}

In this section, we exhibit an asymptotic expansion for $\uu_{\eps}$ and  $\sp_{\eps}$, \S \ref{SFAE}.  
%elementary problems satisfied by  
%formal calculations to derive an asymptotic expansion 
We explicit the first terms in asymptotics, \S \ref{ft}. Then we construct formally equivalent conditions, \S \ref{CECs}.  In \S \ref{SEER} we validate the asymptotic expansion with estimates for the remainders. The main result of this section is the Theorem  \ref{P1} in \S \ref{VECs} which proves the stability of equivalent conditions.

\subsection{Multiscale expansion}
\label{SFAE}
We can exhibit series expansions in powers of $\eps$ for $\uu_{\eps}$ and $\sp_{\eps}$ :  
\begin{gather}
\label{Euu}
   \uu_{\eps}(\xx) \approx \sum_{j\geqslant0} \eps^j \uu_j(\xx) 
%   = \uu_0(\xx) + \eps\uu_1(\xx) + \eps^2 \uu_2(\xx)+ \cdots
   %\mathcal{O}(\eps^3)
   \, ,\\
   \label{Esp}
   \sp_{\eps}(\xx)  \approx \sum_{j\geqslant0} \eps^j   \sp_j(\xx;\eps)
   \quad\mbox{with}\quad       \sp_j(\xx;\eps) = \fp_j(y_{\alpha},\frac{y_{3}}{\eps})\, ,
\end{gather}
%In \eqref{Euu}-\eqref{Esp}, the symbol $\mathcal{O}(\eps^2)$ means that the remainder is uniformly bounded by $\eps^2$.
see Sec. \ref{SEER} for precise estimates. Here $ (y_{\alpha}, y_{3})$ is a ``{\em normal coordinate system}'' ~\cite {CDFP11,Pe13} to the surface $\Gamma$ on the manifold $\Omega^{\eps}_\mathsf{f}$ : $y_{\alpha}$ ($\alpha\in\{1,2\}$) is a tangential coordinate on $\Gamma$ and $y_{3} \in(0,\eps)$ is the distance to the surface $\Gamma$.
The term $\fp_j$ is a ``{\em profile}'' defined on $\Gamma\times (0,1)$. The formal calculus concerning the problem are presented in Sec. \ref{fc} and the first terms $(\fp_j,\uu_j)$ for $j=0,1,2,3$ are explicited in Sec. \ref{ft}.

\subsubsection*{Expansion of the Helmholtz operator}
%\label{not}
%\Rd
%We parameterized the smooth interface $\Gamma$ with the curvilinear abscissa $t\in [0,L)$, such that $\Gamma=\{\xx =\xx(t) \in \R^2 \ |\quad t\in [0,L)\}$. Then, $\exists \eps_{0}>0$ such that $\forall \eps\leqslant \eps_{0}$ the thin layer $\Omega^{\eps}_\mathsf{f}$ can be  parameterized  with curvilinear coordinates $(t,s)$ :
%$$
%\Omega^{\eps}_\mathsf{f}=\{ \xx = \xx(t) + s \nn (t)\in \R^2\ |\quad  t \in \Gamma, \quad s \in (0, \eps)\}\ ,
%$$
%where $\nn (t)=\nn\left(\xx(t)\right)$ denotes the normal vector on $\Gamma$ at the point $\xx(t)$. \Bk
It is possible to write the three dimensional Helmholtz operator in the layer $\Omega^{\eps}_\mathsf{f}$ through the local coordinates $(y_{\alpha}, y_{3})$%there holds for all $y_{3}=h\in (0, \eps)$
%\begin{equation}
%\Delta+ \kappa^2 \mathrm{Id} =(\partial^h_{3})^2  - b_{\alpha}^\alpha (h) \partial^h_{3}  + a^{\alpha\beta}(h) D^h_{\alpha} \partial_{\beta} + \kappa^2 \mathrm{Id} \quad \mbox{on} \quad  \Gamma_{h}\ ,
%\end{equation}
~\cite[Prop.\, B.1]{Pe13}. %Here $\partial_{3}^h$ is  the partial derivative with respect to the normal coordinate $h$, $D^h_{\alpha}$ is the covariant derivative on the manifold $\Gamma_{h}$, $a^{\alpha \beta}(h)$ is the inverse of the  metric tensor in $\Gamma_{h}$, $ b_{\alpha}^{\alpha}(h)= a^{\alpha \beta}(h) b_{\beta\alpha}(h)$ and $b_{\alpha\beta}(h)$ is the {\em curvature tensor} in $\Gamma_{h}$. 
 Then we make the scaling $Y_{3}=\eps^{-1} y_{3} \in(0,1)$ into the normal coordinate  %It maps $\Gamma\times (0,\eps)$ onto  $\Gamma\times (0,1)$. 
% The parameter $\eps$ does not appear anymore in the geometry of the problem but in the equations \eqref{EA}.
 and we expand formally the Helmholtz operator in power series of $\eps$ with coefficient intrinsic operators :
$$
 \Delta + \kappa^2 \mathrm{Id} = \eps^{-2}\left( \di\sum_{n= 0}^{N-1}\eps^n \L^{n}  
+\eps^{N} \mathrm{R}^N_{\eps} \right) \quad \mbox{for all} \quad N\in\N^* \ .
$$
The operators $\L^{n}$, $n=0,1,2$, are explicited in~\cite[Prop.\, B.3]{Pe13}  :
$$
\L^0= \partial_{3}^2 \ , \quad  \L^1= 2 \cH(y_{\alpha})  \partial_{3}\ ,  \quad \mbox{and}\quad \L^2 =\Delta_{\Gamma}+ \kappa^2 \Id- 2\left(2 \cH^2-\cK\right)(y_{\alpha})  Y_{3} \partial_{3}\ .
$$
 Here $\partial_3$ is the partial derivative with respect to $Y_{3}$. We remind that $\Delta_{\Gamma}$ is the Laplace-Beltrami operator along $\Gamma$, $\cH$  
 %=-\frac12 b_{\alpha}^\alpha$  \cite[\S 3.1]{CDFP11}
  and $\cK$ are the {\em mean} and the {\em Gaussian curvature} of the surface $\Gamma$. 
  %The generic operator $\L^{n}$ has smooth coefficients in $y_{\alpha}$ and polynomial in $Y_{3}$ of degree $n-1$; it contains at most one differentiation with respect to $Y_{3}$. 
  The remainder $ \mathrm{R}_{\eps}^N $  has smooth coefficients in $y_{\alpha}$ and $Y_{3}$ which are bounded in $\eps$. 
%  \Bl
% \begin{rem}
% The sign of $\cH$ depends on the orientation of the surface $\Gamma$. As a convention, the unit normal vector $\nn$ on $\Gamma$ is outwardly oriented to $\Omega_{\s}$, see Figure~\ref{FCylGeo}.
 %\end{rem} 
 %\Bk
% \Bl
%In the thin layer, the normal derivative writes 
%$$
 %\partial_{\nn }=\eps^{-1} \partial_{3}\ .
 %$$
%For any function $\sp$ defined in $ \Omega^{\eps}_\mathsf{f}$, we denote  by $\fp$ the function such that :
%$$
 %\sp(\xx) = \fp(y_{\alpha},Y_{3}) \ , \quad (y_{\alpha},Y_{3}=\frac{h}{\eps})\in \Gamma \times (0,1) \ .
%$$
%\Bk

\subsubsection{Elementary problems}
\label{fc}
After the change of variables $y_{3} \mapsto Y_{3}=\eps^{-1} y_{3}$ in the thin layer $\Omega^{\eps}_\mathsf{f}$, problem \eqref{EA} becomes :

\begin{equation}
\label{EA_S}
 \left\{
   \begin{array}{lll}
   \eps^{-2} [\partial_{3}^2  \fp_{\eps} + \di\sum_{n\geqslant 1}\eps^n \L^{n}  \fp_\eps ]
= 0   \quad&\mbox{in}\quad \Gamma\times (0,1)
\\[0.5ex]
\eps^{-1}\partial_{3}\fp_{\eps}=\rho_{\sf} \omega^2 \uu_{\eps}\cdot \nn \quad&\mbox{on}\quad \Gamma\times \{0\}
\\[0.2ex]
\fp_{\eps}= 0\quad&\mbox{on}\quad \Gamma \times \{1\}
% \partial_{\nn }\sp_{\eps} - i \kappa\sp_{\eps}=0 or Dirichlet
\\[0.5ex]
\nabla \cdot \underline{\underline{\sigma}}( \uu_{\eps}) +\omega^2 \rho \uu_{\eps} =\ff
 \quad&\mbox{in}\quad \Omega_{\mathsf{s}}
\\[0.5ex]
\TT(\uu_{\eps}) = -\fp_{\eps} \nn   \quad &\mbox{on} \quad \Gamma \ .
   \end{array}
    \right.
\end{equation}

Inserting the Ansatz \eqref{Euu}-\eqref{Esp} in equations \eqref{EA_S}, we get the following two families of problems, coupled by their boundary conditions on $\Gamma$ (i.e. when $Y_{3}=0$):

\begin{equation}
\label{EAn+}
 \left\{
   \begin{array}{lll}
   \partial_{3}^2  \fp_{n} =-  \di\sum_{l+m=n,l\geqslant 1} \L^{l}  \fp_{m}
    \quad&\mbox{for}\quad Y_{3}\in (0,1)
\\[0.5ex]
 \partial_{3}\fp_{n}=\rho_{\sf} \omega^2 \uu_{n-1}\cdot \nn \quad&\mbox{for} \quad Y_{3}=0
\\[0.2ex]
\fp_{n}= 0\quad&\mbox{for} \quad Y_{3}=1
  \end{array}
    \right.
\end{equation}
\begin{equation}
\label{EAn-}
 \left\{
   \begin{array}{lll}
\nabla \cdot  \underline{\underline{\sigma}}( \uu_{n})+\omega^2 \rho \uu_{n} =\ff \delta_{0}^n
 \quad&\mbox{in}\quad \Omega_{\mathsf{s}}
\\[0.5ex]
\TT(\uu_{n}) = -\fp_{n} \nn   \quad &\mbox{on} \quad \Gamma \ .
   \end{array}
    \right.
\end{equation}
In \eqref{EAn+}, we use the convention $\uu_{-1}=0$, and in \eqref{EAn-}   $\delta_{0}^n$ denotes the Kronecker symbol.

\subsection{First terms}
\label{ft}
We find successively from \eqref{EAn+}-\eqref{EAn-} when $n=0,1,2,3$
%In the case $n=0$, we obtain from \eqref{EAn+} 
$$
\fp_{0}=0 \ ,
$$

%and then \eqref{EAn-} yields $\uu_{0}$ solves the problem%the elastic equation in $  \Omega_{\mathsf{s}}$ with $\TT(\uu_{0}) =-\sp_{i} \nn$:
\begin{equation}
\label{EA0-}
 \left\{
   \begin{array}{lll}
\nabla \cdot \underline{\underline{\sigma}}( \uu_{0}) +\omega^2 \rho \uu_{0} =\ff
 \quad&\mbox{in}\quad \Omega_{\mathsf{s}}
\\[0.5ex]
\TT(\uu_{0}) = 0 \quad &\mbox{on} \quad \Gamma \ ,
   \end{array}
    \right.
\end{equation}

 %At step $n=1$, we find successively 
 $$
 \fp_{1}(., Y_3)=(Y_{3}-1) \rho_{\sf} \omega^2 \uu_{0}\cdot \nn |_\Gamma\ ,
 $$ 
 
 %and that $\uu_{1}$ solves %the elastic equation in $  \Omega_{\mathsf{s}}$ with $\TT(\uu_{1}) = \rho_{\sf} \omega^2 \uu_{0}\cdot \nn\ \nn- \partial_{\nn }\sp_{i}\nn$:
 \begin{equation}
\label{EA1-}
 \left\{
   \begin{array}{lll}
\nabla \cdot  \underline{\underline{\sigma}}(\uu_{1})+\omega^2 \rho \uu_{1} =0
 \quad&\mbox{in}\quad \Omega_{\mathsf{s}}
\\[0.5ex]
\TT(\uu_{1}) = \rho_{\sf} \omega^2 \uu_{0}\cdot \nn\ \nn \quad &\mbox{on} \quad \Gamma \ ,
   \end{array}
    \right.
\end{equation}

%At step $n=2$, since $\fp_{0}=0$, we find 
%       $$
 \begin{equation}
\label{EA2+}
 \fp_{2}(., Y_3)= - (Y_{3}^2-1) \cH \rho_{\sf} \omega^2 \uu_{0}\cdot \nn|_\Gamma
+(Y_{3}-1)\rho_{\sf} \omega^2 \uu_{1}\cdot \nn|_\Gamma \ ,
\end{equation}

%and then, $\uu_{2}$ solves the problem
 \begin{equation}
\label{EA2-}
 \left\{
   \begin{array}{lll}
\nabla \cdot \underline{\underline{\sigma}}(\uu_{2})+\omega^2 \rho \uu_{2} =0
 \quad&\mbox{in}\quad \Omega_{\mathsf{s}}
\\[0.5ex]
\TT(\uu_{2})= \rho_{\sf} \omega^2 \left( \uu_{1}
-\cH   \uu_{0}\right) \cdot \nn\ \nn 
 \quad &\mbox{on} \quad \Gamma \ ,
   \end{array}
    \right.
\end{equation}

%  At step $n=3$, we find 
 \begin{multline}
\label{EA3+}
 \fp_{3}(., Y_3)= (Y_{3}^3-1)\frac16 \rho_{\sf} \omega^2 \left( 8\cH^2-2\cK - (\Delta_{\Gamma} +\kappa^2 \Id) \right) (\uu_{0}\cdot \nn)
  \\
 +  (Y_{3}^2-1)\frac12  \rho_{\sf} \omega^2  \left( (\Delta_{\Gamma} +\kappa^2 \Id) (\uu_{0}\cdot \nn) - 2\cH  \uu_{1}\cdot \nn\right)
 + (Y_{3}-1) \rho_{\sf}\omega^2 \uu_{2}\cdot \nn \ ,
 \end{multline}
 
%Then, $\uu_{3}$ solves the problem
 \begin{equation}
\label{EA3-}
 \left\{
   \begin{array}{lll}
\nabla \cdot  \underline{\underline{\sigma}}(\uu_{3})+\omega^2 \rho \uu_{3} =0
 \quad\mbox{in}\quad \Omega_{\mathsf{s}}
\\[0.5ex]
\TT(\uu_{3})= \rho_{\sf} \omega^2 \left( \frac13\left( 4\cH^2- \cK +\Delta_{\Gamma} +\kappa^2 \Id \right) (\uu_{0}\cdot \nn)+
( -\cH \uu_{1}+
  \uu_{2})\cdot \nn
  \right) \nn  \ \ \ \mbox{on}\ \ \Gamma\ .
   \end{array}
    \right.
\end{equation}

We refer the reader to~\cite{Pe13} for more details. The whole construction of the asymptotics comes from a recurrence argument :
%, see for instance ~\cite {CCDV06} for a similar context. 
if the sequences $(\uu_{n})$ and $(\fp_{n})$ are known until rank $n=N-1$, then the Sturm-Liouville problem  \eqref{EAn+} uniquely defines $\fp_{N}$ whose trace on $\Gamma$ is inserted into  \eqref{EAn-}  as a data to determine $\uu_{N}$. The next proposition ensures
% elementary problems are well-posed together with optimal
 regularity results for the first terms $\uu_k$ in $\bV^{k}$ (Not. \ref{not2}) $k=0,1,2,3$. %Remind that $\bV^{k}=\bH^1(\Omega_{\mathsf{s}})$ when $k\in\{0,1,2\}$ and  $\bV^{3}=\{\uu\in \bH^1(\Omega_{\mathsf{s}}) \, | \quad \uu\cdot \nn|_{\Gamma} \in\H^1(\Gamma) \}$.

\begin{prop}
Under Assumptions \ref{AC}-\ref{SA}-\ref{HG}, if $\ff\in\bL^2(\Omega_{\s})$ then  elementary problems \eqref{EA0-}-\eqref{EA1-}-\eqref{EA2-}-\eqref{EA3-} have a unique solution $\uu_k$ in $\bV^k$ when $k=0,1,2,3$. 

Let $l$ be a non-negative integer. If $\Omega_{\s}$ is of class $\mathcal{C}^{l+2}$ and $\ff\in\bH^{l}(\Omega_{\s})$, then $\uu_{0}$, $\uu_{3}$ belong to $\bH^{l+2}(\Omega_{\s})$ and $\uu_{1}$, $\uu_{2}$ belong to $\bH^{l+3}(\Omega_{\s})$.
\end{prop}
Here, the regularity of $\uu_{0}$, $\uu_{1}$, $\uu_{2}$ and $\uu_{3}$ is a consequence of a general shift result available in Sobolev spaces~\cite[Th. 3.4.5]{CDN10}, since $\fp_{1}|_{\Gamma}\in \H^{l+\frac32}(\Gamma)$, $\fp_{2}|_{\Gamma}\in \H^{l+\frac32}(\Gamma)$ and $\fp_{3}|_{\Gamma}\in \H^{l+\frac12}(\Gamma)$.

\subsection{Construction of equivalent conditions}
\label{CECs}
In this section, we derive formally ECs (Sec. \ref{S3DECs}).
\subsubsection*{Order 1}
Since the equations in \eqref{EA0-} are independent of $\eps$, the condition of order $1$ is  the stress free boundary condition
$$
\TT(\uu_{0})= 0 \quad \mbox{on} \quad \Gamma\ .
$$
%This condition  is "exact" i.e. $\uu_\eps^{0}=\uu_{0}$. 
\subsubsection*{Order 2}
According to \eqref{EA0-} and \eqref{EA1-}, the truncated expansion $\uu_{1,\eps}=\uu_{0}+\eps\uu_{1}$  solves the elastic equation in $  \Omega_{\mathsf{s}}$ together with the boundary condition 
$$
\TT(\uu_{1,\eps})= \eps \rho_{\sf} \omega^2 \uu_{0}\cdot \nn\ \nn  \quad \mbox{on} \quad \Gamma\ .
$$
Writting $ \uu_{0}= \uu_{1,\eps}-\eps \uu_{1}$, there holds
%$$
%\TT(\uu_{1,\eps})=-\sp_{i} \nn+\eps (\rho_{\sf} \omega^2( \uu_{1,\eps}-\eps \uu_{1})\cdot \nn\ \nn- \partial_{\nn }\sp_{i}\nn) \quad \mbox{on} \quad \Gamma\ ,
%$$
%i.e.
$$
\TT(\uu_{1,\eps}) -\eps \omega^2 \rho_{\sf} \uu_{1,\eps}\cdot \nn \, \nn= -\eps^2 \rho_{\sf} \omega^2 \uu_{1}\cdot \nn\ \nn   \quad \mbox{on} \quad \Gamma\ .
$$
Neglecting the term of order $\eps^2$ in the previous right-hand side, we infer the  condition \eqref{Ebc1}.
%we obtain the asymptotic model of order $2$: we define $\uu_\eps^{1}$ such that it solves the elastic equation in $  \Omega_{\mathsf{s}}$ together with the equivalent condition
%$$
%\TT(\uu_\eps^{1}) -\eps \omega^2 \rho_{\sf} \uu_\eps^{1}\cdot \nn \, \nn= 0 \quad \mbox{on} \quad \Gamma\ .
%$$
%\begin{rem}
%$\uu_\eps^{1}-\uu_{1,\eps}$ solves the elastic equation in $  \Omega_{\mathsf{s}}$ with the boundary condition
%$$
%\TT(\uu_\eps^{1}-\uu_{1,\eps}) -\eps \omega^2 \rho_{\sf}(\uu_\eps^{1}-\uu_{1,\eps}) \cdot \nn \, \nn=  \eps^2 \rho_{\sf} \omega^2 \uu_{1}\cdot \nn\ \nn \quad \mbox{on} \quad \Gamma\ .
%$$
%\end{rem}

\subsubsection*{Order 3}
According to \eqref{EA0-}, \eqref{EA1-}, and \eqref{EA2-}, the truncated expansion $\uu_{2,\eps}=\uu_{0}+\eps\uu_{1}+\eps^2\uu_{2}$  solves the elastic equation in $  \Omega_{\mathsf{s}}$ with the condition %set on $\Gamma$
%$$
%\TT(\uu_{2,\eps})=-\sp_{i} \nn+\eps (\rho_{\sf} \omega^2 \uu_{0}\cdot \nn\ \nn- \partial_{\nn }\sp_{i}\nn)
%+\eps^2\left(
 %\rho_{\sf} \omega^2 \uu_{1}\cdot \nn \ \nn
%+\cH [\rho_{\sf} \omega^2 \uu_{0}\cdot \nn- \partial_{\nn }\sp_{i}] \nn 
%\right)
% \quad \mbox{on} \quad \Gamma
%\ ,
%$$ 
%i.e.
%\begin{equation}
%\label{Tu2eps}
%\TT(\uu_{2,\eps})=\eps \rho_{\sf} \omega^2 \left (\uu_{0}+\eps \uu_{1}\right) \cdot \nn\ \nn
%- \eps^2 \cH   \rho_{\sf} \omega^2 \uu_{0}\cdot \nn \ \nn 
% \quad \mbox{on} \quad \Gamma\ .
%\end{equation}
%Using $ \uu_{0}+\eps \uu_{1}= \uu_{2,\eps}-\eps^2 \uu_{2}$ and $ \uu_{0}= \uu_{2,\eps}-\eps \uu_{1}- \eps^2 \uu_{2}$, we obtain on $\Gamma$:
$$
\TT(\uu_{2,\eps})- \eps \omega^2 \rho_{\sf} 
\left(1-\eps \cH  \right)  \uu_{2,\eps}\cdot \nn \, \nn =  - \eps^3 \omega^2 \rho_{\sf} \uu_{2}\cdot \nn \, \nn  +\eps^3 \cH  \omega^2 \rho_{\sf}  (\uu_{1}+ \eps \uu_{2})\cdot \nn \, \nn  \ .
%\quad \mbox{on} \quad \Gamma
$$
%Neglecting the terms of order $\eps^3$ in the right-hand sides of the previous condition, we obtain the equivalent condition of order $3$ %$\uu_\eps^{2}$ solves the elastic equation in $  \Omega_{\mathsf{s}}$ with the equivalent condition
Neglecting the term of order $\eps^3$ in the right-hand side, we obtain the  condition \eqref{Ebc2}.
% $$
%\TT(\uu_\eps^{2})- \eps \omega^2 \rho_{\sf} 
%\left(1-\eps \cH  \right) \uu_\eps^{2}\cdot \nn \, \nn = 0 \quad \mbox{on} \quad \Gamma\ .
%$$
%\begin{rem}
%$\uu_\eps^{2}-\uu_{2,\eps}$ solves the elastic equation in $  \Omega_{\mathsf{s}}$ with the b.c. 
%\begin{multline*}
%\TT(\uu_\eps^{2}-\uu_{2,\eps}) -\eps \omega^2 \rho_{\sf}(\uu_\eps^{2}-\uu_{2,\eps}) \cdot \nn \, \nn
%+\frac{\eps^2}{2} \mathrm{c}(t) \omega^2 \rho_{\sf}  (\uu_\eps^{2}-\uu_{2,\eps})\cdot \nn \, \nn
%=  \eps^3 \omega^2 \rho_{\sf} \uu_{2}\cdot \nn \, \nn  
%\\
%-\frac{\eps^3}{2} \mathrm{c}(t) \omega^2 \rho_{\sf}  (\uu_{1}+ \eps \uu_{2})\cdot \nn \, \nn   \quad \mbox{on} \quad \Gamma\ .
%\end{multline*}
%\end{rem}

\subsubsection*{Order 4}
According to \eqref{EA0-}, \eqref{EA1-}, \eqref{EA2-} and \eqref{EA3-}, the truncated  expansion $\uu_{3,\eps}=\uu_{0}+\eps\uu_{1}+\eps^2\uu_{2}+\eps^3\uu_{3}$  solves the elastic equation in $  \Omega_{\mathsf{s}}$ with the condition %set on $\Gamma$
%\begin{equation}
%\TT(\uu_{3,\eps})= \TT(\uu_{2,\eps})+\eps^3 \rho_{\sf} \omega^2  \large[
%\frac13\left( 4\cH^2-\cK +\Delta_{\Gamma} +\kappa^2 \Id \right) \uu_{0}
%+ (- \cH \uu_{1}+\uu_{2}) \large ] \cdot \nn\ \nn
 %\ .
 %\end{equation}
%According to \eqref{Tu2eps}, there holds
%\begin{multline*}
%\TT(\uu_{3,\eps})=\eps \rho_{\sf} \omega^2 \left (\uu_{0}+\eps \uu_{1}+\eps^2 \uu_{2}\right) \cdot \nn\ \nn
%+ \eps^2 \cH   \rho_{\sf} \omega^2 (\uu_{0}+\eps \uu_{1})\cdot \nn\ \nn 
%\\
%+\eps^3 
%\frac16\left( 8\cH^2-\cK+2(\Delta_{\Gamma} +\kappa^2 \Id) \right) \rho_{\sf} \omega^2 \uu_{0}\cdot \nn\ \nn
% \quad \mbox{on} \quad \Gamma\ .
%\end{multline*}
%We infer
\begin{multline*}
\TT(\uu_{3,\eps}) - \eps \omega^2 \rho_{\sf} 
\left[1-\eps \cH  +\frac{\eps^2}{3}\left( 4\cH^2-\cK+\Delta_{\Gamma} +\kappa^2 \Id  \right)
\right] ( \uu_{3,\eps}\cdot \nn) \, \nn=
%\hh_{3,\eps}
\\
-\eps^4 \rho_{\sf} \omega^2\left[
\left( \uu_{3}- \cH  (\uu_{2}+ \eps \uu_{3})\right)\cdot \nn
+\frac13\left( 4\cH^2-\cK + \Delta_{\Gamma} +\kappa^2 \Id \right)
 \left((\uu_{1}+ \eps \uu_{2}+\eps^2 \uu_{3})\cdot \nn\right)
\right]   \nn  \ .
%-\sp_{i} \nn
%+\eps (\rho_{\sf} \omega^2 \left (\uu_{0}+\eps \uu_{1}+\eps^2 \uu_{2}\right) \cdot \nn\ \nn- \partial_{\nn }\sp_{i}\nn)
%+ \eps^2 \cH   [\rho_{\sf} \omega^2 (\uu_{0}+\eps \uu_{1})\cdot \nn- \partial_{\nn }\sp_{i}] \nn 
%\\
%+\eps^3 [
%\frac13\left( \frac12(8\cH^2-\cK) +\Delta_{\Gamma} +\kappa^2 \Id \right) (\rho_{\sf} \omega^2 \uu_{0}\cdot \nn- \partial_{\nn }\sp_{i})
 %\quad \mbox{on} \quad \Gamma\ .
\end{multline*}
%with
%$$
%\hh_{3,\eps}= -\sp_{i} \nn -\eps \left[1 + \eps \cH + 
%\eps^2 \frac16\left( 8\cH^2-\cK +2 (\Delta_{\Gamma} +\kappa^2 \Id) \right)
%\right] (\partial_{\nn }\sp_{i}) \nn\ .
%$$
 %Neglecting the terms of order $\eps^4$ in the right-hand sides of the previous condition, we obtain the equivalent condition of order $4$%: $\uu_\eps^{3}$ solves the elastic equation in $  \Omega_{\mathsf{s}}$ with the equivalent condition
 Neglecting the previous right-hand side, we infer the  condition \eqref{Ebc3}.
% $$
%\TT(\uu_\eps^{3})
% - \eps \omega^2 \rho_{\sf} 
%\left[1-\eps \cH  +\frac{\eps^2}{3}\left( 4\cH^2-\cK+\Delta_{\Gamma} +\kappa^2 \Id \right)
%\right] ( \uu_{\eps}^3\cdot \nn) \, \nn= 0
%\hh_{3,\eps}
% \quad \mbox{on} \quad \Gamma\ .
%$$

\subsection{Estimates for the remainders}
\label{SEER}

The validation of the asymptotic expansion  \eqref{Euu}-\eqref{Esp} consists in proving estimates for remainders $(\mathbf{ r}_{\eps}^{N}, r_{\eps}^{N})$ defined in $\Omega_{\mathsf{s}}$ and $\Omega^{\eps}_\mathsf{f}$ as 

%By construction, the remainder $(\mathbf{ r}_{\eps}^{N}, r_{\eps}^{N}):= (\uu_{\eps},\sp_{\eps})-\di\Sigma_{n=0}^N \eps^n (\uu_{n},\sp_{n})$ defined as
\begin{equation}
\label{remN}
\mathbf{ r}_{\eps}^{N}= \uu_{\eps}-\sum_{n=0}^N \eps^n \uu_{n} \quad\mbox{in}\quad \Omega_{\mathsf{s}}\ ,
\quad\mbox{and} \quad
 r_{\eps}^{N}(\xx)= \sp_{\eps}(\xx)-\sum_{n=0}^N \eps^n  \fp_n(y_{\alpha},\frac{y_{3}}{\eps})
 %\sp_{n}(\xx;\eps) 
 \quad\mbox{for all}\, \xx\in\Omega^{\eps}_\mathsf{f}\ .
\end{equation}
The convergence result is the following statement.
\begin{thm}
\label{TErem}
%Let $\sp_{i}$ belong to $ \mathcal{C}^\infty(\Gamma)$.
Under Assumptions \ref{AC}-\ref{SA}-\ref{HG}-\ref{Hspi} and for $\eps$ small enough, the solution $(\uu_{\eps},\sp_{\eps})$ of problem \eqref{EA} 
%with Dirichlet external b.c.
 has  a two-scale expansion  which can be written  in the form  \eqref{Euu}-\eqref{Esp}, with $\uu_{j}\in\bH^1(\Omega_{\s})$ and $\fp_{j}\in\H^1\left(\Gamma\times(0,1)\right)$. For each $N\in\N$, the remainders $(\mathbf{ r}_{\eps}^{N}, r_{\eps}^{N})$ satisfy %optimal estimates
\begin{equation}
\label{ErN}
 \|  \mathbf{ r}_{\eps}^{N} \|_{1,\Omega_{\s}} 
  +{{\sqrt{\eps}}}\| r_{\eps}^{N}\|_{1,{\Omega^{\eps}_\mathsf{f}}}  \leqslant C_{N} {{\eps^{N+1}}}
\end{equation}
with a constant $C_{N}$ independent of $\eps$.
\end{thm}

%The proof of this Theorem is detailed in~\cite{Pe13}. 
The error estimate \eqref{ErN} is obtained through an evaluation of the right-hand sides when applying Theorem \ref{tUEeps} to the couple $(\uu,\sp)=(\mathbf{ r}_{\eps}^{N}, r_{\eps}^{N})$.

\begin{proof}
The proof is rather standard, see for instance the proof of \ \cite[Th.\,2.1]{CCDV06} where the authors consider an interface problem for the Laplacian operator set in a domain with a thin layer. The error estimate \eqref{ErN} is obtained through an evaluation of the right-hand sides when the elasto-acoustic operator is applied to $(\mathbf{ r}_{\eps}^{N}, r_{\eps}^{N})$. By construction, the remainder $(\mathbf{ r}_{\eps}^{N}, r_{\eps}^{N})$
 is solution of problem
 \begin{equation}
\label{Erem}
 \left\{
   \begin{array}{lll}
 \Delta r_{\eps}^{N} + \kappa^2 r_{\eps}^{N}= 
 f_{N,\eps} \quad&\mbox{in}\quad \Omega^{\eps}_\mathsf{f}
\\[0.8ex]
 \nabla \cdot \underline{\underline{\sigma}}( \mathbf{ r}_{\eps}^{N} )+\omega^2 \rho  \mathbf{ r}_{\eps}^{N} =0
 \quad&\mbox{in}\quad \Omega_{\mathsf{s}}
\\[0.8ex]
\partial_{\nn } r_{\eps}^{N} =\rho_{\sf} \omega^2  \mathbf{ r}_{\eps}^{N}  \cdot \nn+ g_{N,\eps}
 \quad&\mbox{on}\quad \Gamma
\\[0.8ex]
\TT(\mathbf{ r}_{\eps}^{N} ) = -r_{\eps}^{N}   \nn \quad &\mbox{on} \quad \Gamma
\\[0.8ex]
%\mbox{for Dirichlet external B.C.}\quad 
 r_{\eps}^{N} =0 
  \quad &\mbox{on}\quad \Gamma^\eps \ . % \partial_{\nn }\sp_{\eps} - i \kappa\sp_{\eps}=0 or Dirichlet
   \end{array}
    \right.
\end{equation}
Here, the right-hand sides are explicit : 
$$
f_{N,\eps}=\eps^{N-1} [ F_{N} -\di\sum_{l=1}^{N} \eps^{l-2+N}  R_{\eps}^{N-l+1} \fp_{l} ]  \quad\mbox{in}\quad \Omega^{\eps}_\mathsf{f}\ ,
$$ 
and 
$$
g_{N,\eps}=
%- \partial_{\nn }\sp_{i} \delta_{0}^N +
\rho_{\sf} \omega^2 \eps^N \uu^{N}  \cdot \nn \quad \mbox{on}\quad \Gamma\ .
$$ 
%\begin{rem}
%For $N=0$,   $\kappa^2 r_{\eps}^{N+1}+ \Delta r_{\eps}^{N+1} =0$. $F_{1}= 0$. $F_{2}= -A_{1} \fp_{2}$.
%\end{rem}
We have the following estimates for the residues $f_{N,\eps}$ and $g_{N,\eps}$ 
\begin{equation*}
\| f_{N,\eps} \|_{0, \Omega^{\eps}_\mathsf{f}}=\mathcal{O}(\eps^{N-\frac12}) \quad \mbox{and} \quad 
\|g_{N,\eps}\|_{0,\Gamma}=\mathcal{O}(\eps^{N})\ .
\end{equation*}
We can apply Theorem \ref{tUEeps} to the couple $(\uu,\sp)=(\mathbf{ r}_{\eps}^{N}, r_{\eps}^{N})$, and we obtain 
$$
 \|  \mathbf{ r}_{\eps}^{N} \|_{1,\Omega_{\s}} 
  +\| r_{\eps}^{N}\|_{1,{\Omega^{\eps}_\mathsf{f}}}  \leqslant C_{N} {{\eps^{N-\frac12}}}\ .
$$
Writting 
$\mathbf{ r}_{\eps}^{N}= \mathbf{ r}_{\eps}^{N+2}+\eps^{N+2} \uu_{N+2}  +\eps^{N+1} \uu_{N+1}$ and  $r_{\eps}^{N}=  r_{\eps}^{N+2}+\eps^{N+2} \fp_{N+2}  +\eps^{N+1} \fp_{N+1}$, we apply the previous estimate to the couple $(\mathbf{ r}_{\eps}^{N+2}, r_{\eps}^{N+2})$
and we use estimates 
$$ 
\| \uu_{l}\|_{1,\Omega_{\s}}=\mathcal{O}(1) \quad \mbox{and} \quad  
\| \fp_{l}\|_{1,\Omega^\eps_{\f}}=\mathcal{O}(\eps^{-\frac12})
$$
 to infer the optimal estimates \eqref{ErN}.
%We conclude as in the proof of \ \cite[Th.\,2.1]{CCDV06}.
\end{proof}

\subsection{Validation of equivalent conditions}
\label{VECs}
%Our interest lies in the elastic wave problem with an equivalent condition set on the boundary of the domain \eqref{Ebck} in time-harmonic regime.

%(here $\uu_\eps=\uu_\eps^{1}$ and $\h \nn=\hh_{1,\eps}=-\sp_{i} \nn   -\eps  \partial_{\nn }\sp_{i} \nn$).

%\Rd
%We assume : we can choose a standard lifting in $\HH^2(\Omega_{\mathsf{s}})$ of the two traces :
%$$
%\exists \vv\in \HH^2(\Omega_{\mathsf{s}}) \quad \mbox{such that} \quad \vv\cdot\nn =0 \quad \mbox{on} \quad \Gamma \ ,\quad \mbox{and} \quad \TT(\vv) =\hh \quad \mbox{on} \quad \Gamma
%$$
%and $ \|\vv \|_{2,\Omega_{\mathsf{s}}}\leqslant \| \hh\|_{0,\Gamma}$\ .
%Then  $\uu:=\uu_\eps-\vv$ satisfies:
%\begin{equation}
%\label{Eec1bis}
 %\left\{
 %  \begin{array}{lll}
%\omega^2 \rho \uu + \nabla \cdot  \underline{\underline{C}} \ \underline{\underline{\nabla}} \uu=\ff
% \quad&\mbox{in}\quad \Omega_{\mathsf{s}}
%\\[0.5ex]
%\TT(\uu) -\eps \omega^2 \rho_{\sf} \uu\cdot \nn \, \nn= 0  \quad &\mbox{on} \quad \Gamma\ ,  
%   \end{array}
 %   \right.
%\end{equation}
%with $\ff=-\omega^2 \rho \vv - \nabla \cdot  \underline{\underline{C}} \ \underline{\underline{\nabla}} \vv$.
%\Bk

We consider the problem \eqref{Ebck} with an equivalent condition and at a fixed frequency $\omega$  satisfying Assumption \ref{SA}. %Recall that some resonant frequencies may appears in the solid domain, hence, we work under the spectral assumption \ref{SA}.
%The functional setting for the problem \eqref{Ebck} is described by the Hilbert space  $\bV^{k}$.
%\begin{itemize}
%\item $\bH^1(\Omega_{\mathsf{s}})$ when $k$ belongs to $\{0,1,2\}$
%\item $\bV=\{\uu\in \bH^1(\Omega_{\mathsf{s}}) \, | \quad \uu\cdot \nn|_{\Gamma} \in\H^1(\Gamma) \}$  for the Ventcel's problem (when $k=3)$ .
%\end{itemize}
%The space $\bV$ endowed with it's natural norm is a Hilbert space. To unify the notations, we introduce a common variational framework  %denoting by
%\begin{itemize}
%\item $\bV^{k}$ denotes the space $\bH^1(\Omega_{\mathsf{s}})$ when $k\in\{0,1,2\}$
%\item $\bV^{3}$ denotes the space $\bV$.
%\end{itemize}
 % =$\bH^1(\Omega_{\mathsf{s}})$ when $k\in\{0,1,2\}$ and  $\bV^{3}=\{\uu\in \bH^1(\Omega_{\mathsf{s}}) \, | \quad \uu\cdot \nn|_{\Gamma} \in\H^1(\Gamma) \}$.
 The main result of this section is the following statement, that is for all $k \in \{0,1,2,3\}$ the problem  \eqref{Ebck} is well-posed in the space $\bV^{k}$ (Not. \ref{not2}), and its solution satisfies uniform $\bH^1$ estimates.

\begin{thm}
\label{P1}
Under  Assumptions \ref{AC}-\ref{SA}-\ref{HG},  for all $ k \in \{0,1,2,3\}$ there are constants $\eps_{k},C_{k}>0$ such that for all $ \eps \in (0, \eps_{k})$, the problem \eqref{Ebck}  with a data $\ff\in\bL^2(\Omega_{\s})$  has a unique solution $\uu^k_{\eps}\in\bV^{k}$ which satisfies the uniform estimates:
\begin{subequations}
\begin{gather}
\label{UE1a}
\|\uu^k_{\eps} \|_{1,\Omega_{\mathsf{s}}}\leqslant C_{k} \| \ff \|_{0,\Omega_{\s}}\quad \mbox{for all} \quad k \in \{0,1,2\} \
 \,,\\[0.8ex]
\label{UE1b}
\|\uu^3_{\eps} \|_{1,\Omega_{\mathsf{s}}} + \eps^{\frac32} \| \nabla_{\Gamma}(\uu^3_{\eps}\cdot\nn)  \|_{0,\Gamma}  \leqslant C_{3}  \| \ff \|_{0,\Omega_{\s}} \ .
\end{gather}
\end{subequations}
%\begin{equation}
%\label{UE1}
%\|\uu^k_{\eps} \|_{1,\Omega_{\mathsf{s}}}\leqslant C_{k} \| \hh_{k}\|_{0,\Gamma}\ .
%\end{equation}
\end{thm}
The key for the proof of Thm. \ref{P1} is the following Lemma. 
\begin{lem}
\label{L1}
Under Assumptions \ref{AC}-\ref{SA}-\ref{HG}, for all $ k \in \{0,1,2,3\}$ there exists constants  $ \eps_{k},C_{k}>0$  such that for all $\eps\in(0,\eps_{k})$, any solution  $\uu^k_{\eps}\in\bV^{k}$ %$\uu^k_{\eps}\in \bL^2(\Omega_{\mathsf{s}})$ 
of problem \eqref{Ebck} with a data $\ff\in\bL^2(\Omega_{\s})$  satisfies the uniform estimate:
\begin{equation}
\label{UE2}
\|\uu^k_{\eps} \|_{0,\Omega_{\mathsf{s}}}\leqslant C_{k} \| \ff \|_{0,\Omega_{\s}}  \ .
\end{equation}
\end{lem}
\begin{rem}
For $k=0$, the Theorem \ref{P1} and the Lemma \ref{L1} hold for all $\eps>0$.
\end{rem}
The Lemma \ref{L1} is proved in Sec. \ref{SWp}. As a consequence of this Lemma, each solution of the problem \eqref{Ebck} satisfies uniform $\bH^1$-estimates  \eqref{UE1a}-\eqref{UE1b}. Then, the proof of the Thm. \ref{P1} is obtained as a consequence of the Fredholm alternative since the problem \eqref{Ebck} is of Fredholm type. One passes from Lemma \ref{L1} to Theorem  \ref{P1} as from Lemma \ref{L2} to Theorem \ref{UE}.

% This last argument is rather classic, we refer for instance to the proof of \ \cite[Cor.\,4.3]{CDP10}, or \ \cite[Lem.\,5.4]{HJN05} for a similar argument.

\section{Analysis of Equivalent Conditions}
\label{SAEC}

In this section, we first prove the Lemma  \ref{L1}, i.e. uniform $\bL^2$-estimate  \eqref{UE2} for the solution of problem \eqref{Ebck}, Sec. \ref{SWp}. In Sec. \ref{SEe}, we prove that the solution $\uu_\eps^{k}$ of problem \eqref{Ebck}  satisfies uniform $\bH^1$ error estimates \eqref{UEec} and we infer the Theorem \ref{CVec}.
 %  For $k=0$, the  proof of Lemma  \ref{L1} is a consequence of  the Assumption \ref{SA} and \Rd the Fredholm alternative \Bk . 
We focus on the proof  of Lemma  \ref{L1}  for $k=3$ only since the proof when $k\in\{0,1,2\}$ is simpler. Hence we consider the problem (here $\uu=\uu^3_\eps$)
\begin{equation}
\label{E1}
 \left\{
   \begin{array}{lll}
 \nabla \cdot \underline{\underline{\sigma}}(\uu)+\omega^2  \rho \uu =\ff
 \quad&\mbox{in}\quad \Omega_{\mathsf{s}}
\\[0.8ex]
\TT(\uu)+ \B_{3,\eps}(\uu\cdot\nn) \nn 
%\eps \omega^2 \rho_{\sf} \uu_\eps\cdot \nn \, \nn
= 0 \quad &\mbox{on} \quad \Gamma\ .  
   \end{array}
    \right.
\end{equation}

To prepare for the proof, we introduce the variational formulation for $\uu_\eps$. If $\uu\in\bV^3$ is a solution of \eqref{E1}, then it satisfies for all $\vv\in\bV^3$:
\begin{multline}
\label{EFV}
  \int_{\Omega_{\mathsf{s}}}\left( \underline{\underline{\sigma}}(\uu) : \underline{\underline{\epsilon}}  (\bar{\vv})  -\omega^2 \rho \uu  \bar{\vv}  \right)\;\dr\xx  
  \ -\eps\omega^2\rho_{\sf} \int_{\Gamma} \cJ_{\eps} \uu \cdot \nn\bar{\vv} \cdot \nn \,\dr \sigma 
    \\
  +\frac{\eps^3}{3}\omega^2\rho_{\sf} \int_{\Gamma}  \nabla_{\Gamma}( \uu \cdot \nn) \nabla_{\Gamma}(\bar{\vv} \cdot \nn) \,\dr \sigma 
=-  \int_{\Omega_{\mathsf{s}}} \ff \cdot  \bar{\vv} \;\dr\xx   \   ,
\end{multline}
where $\cJ_{\eps}$ is a function defined on $\Gamma$ as $\cJ_{\eps}=\left(1 -\eps \cH +\frac{\eps^2}{3}\left( 4 \cH^2 -\cK -\kappa^2 \right)  \right)$ which tends to  $1$ when $\eps$ goes to  $0$.

\subsection{Proof of Lemma \ref{L1}\,: Uniform $\bL^2$-estimate of the elastic displacement}
\label{SWp}
Reductio ad absurdum: We assume that there is a sequence $(\uu_m)\in\bL^2(\Omega_{\mathsf{s}})$, $m\in\N$, of solutions of the problem \eqref{E1} associated with a parameter $\eps_m$ and a right-hand side $\ff_m\in\bL^2(\Omega_{\s})$:
\begin{subequations}
\begin{gather}
\label{E3a}
 \nabla \cdot  \underline{\underline{\sigma}}(\uu_m)+  \omega^2 \rho \uu_m =\ff_{m}
 \quad\mbox{in}\quad\Omega_{\mathsf{s}}\,,\\
\label{E3b}
  \TT(\uu_{m}) -\eps_{m} \omega^2 \rho_{\sf} \cJ_{\eps_m} \uu_{m}\cdot \nn \, \nn
-  \frac{\eps_{m}^3}{3}\omega^2\rho_{\sf} \Delta_{\Gamma}( \uu_{m} \cdot \nn)\nn
  =  0\quad\mbox{on}\quad\Gamma\,,
\end{gather}
\end{subequations}
satisfying the following conditions
\begin{subequations}
\begin{eqnarray}
\label{5E4a}
   &\eps_m\to 0\quad &\mbox{as \ $m\to\infty$,} \\
\label{5E4b}
   &\|\uu_m\|_{0,\Omega_{\mathsf{s}}} = 1\quad &\mbox{for all $ m\in\N$,} \\
\label{5E4c}
   &\|\ff_m\|_{0,\Omega_{\mathsf{s}}}\to0\quad&\mbox{as \ $m\to\infty$.}
\end{eqnarray}
\end{subequations}

\subsubsection{Estimates of the sequence $\{\uu_{m}\}$}

We first prove that the sequence $\{\uu_m\}$ is bounded (only) in $\bH^1(\Omega_{\mathsf{s}})$
%(but not in $\bV^3$)
. We particularize the variational formulation \eqref{EFV} for the sequence 
$\{\uu_m\}$: For all $\vv\in\bV^3$,
\begin{multline}
\label{EVm}
  \int_{\Omega_{\mathsf{s}}}\left( \underline{\underline{\sigma}}(\uu_{m}) : \underline{\underline{\epsilon}}  (\bar{\vv})  -\omega^2 \rho \uu_{m}  \bar{\vv}  \right)\;\dr\xx  \ -\eps_{m}\omega^2\rho_{\sf} \int_{\Gamma}  \cJ_{\eps_m}  \uu_{m} \cdot \nn\ \bar{\vv} \cdot \nn \,\dr \sigma
\\  
  +\frac{\eps_{m}^3}{3}\omega^2\rho_{\sf} \int_{\Gamma}  \nabla_{\Gamma}( \uu_{m} \cdot \nn) \nabla_{\Gamma}(\bar{\vv} \cdot \nn) \,\dr \sigma 
=-  \int_{\Omega_{\mathsf{s}}} \ff_{m} \cdot  \bar{\vv} \;\dr\xx  \   .
\end{multline}
Choosing $\vv=\uu_m$ in \eqref{EVm}, we obtain with the help of condition \eqref{5E4b} the uniform bound 
\begin{multline}
\label{5EC}
 \int_{\Omega_{\mathsf{s}}}  \underline{\underline{C}} \ \underline{\underline{\epsilon}} (\uu_{m}) : \underline{\underline{\epsilon}} (\bar\uu_{m} )
\;  \dr\xx 
  -\eps_{m}\omega^2\rho_{\sf} \int_{\Gamma}   \cJ_{\eps_m}  |\uu_{m} \cdot \nn|^2 \,\dr \sigma
  \\
  +\frac{\eps_{m}^3}{3}\omega^2\rho_{\sf} \int_{\Gamma}|  \nabla_{\Gamma}( \uu_{m} \cdot \nn) |^2 \,\dr \sigma 
    \leqslant \omega^2 \rho +\|\ff_m\|_{0,\Omega_{\mathsf{s}}} \,.
\end{multline}
%with $\kappa= \omega^2 \rho_{\sf}$. 

Since the tensor $\underline{\underline{\epsilon}}  (\uu )$ is symmetric, thanks to the assumptions \ref{AC} $\iti1$-$\iti3$ together with the Korn inequality%(which is available since $\Omega_{\mathsf{s}}$ is a smooth domain)
, we infer : there exists constants $C,c>0$ such that for all $\uu\in\bH^1(\Omega_{\mathsf{s}})$
\begin{equation}
\label{korn}
 \int_{\Omega_{\mathsf{s}}}  \underline{\underline{C}} \ \underline{\underline{\epsilon}} (\uu) : \underline{\underline{\epsilon}} (\bar\uu )
  \dr\xx  \geqslant \alpha C \| \uu\|^2_{1,\Omega_{\mathsf{s}}} -\alpha c  \| \uu\|^2_{0, \Omega_{\mathsf{s}}}\ .
\end{equation}
%  \begin{rem}
 %  The tensor  $\underline{\underline{\epsilon}}  (\uu )$ is symmetric, hence thanks to the  assumptions \ref{AC} $\iti1$-$\iti3$ there holds
 %      $$ 
  %     \forall \uu\in \bH^1(\Omega_{\mathsf{s}})\ ,\quad \int_{\Omega_{\mathsf{s}}}  \underline{\underline{C}} \ \underline{\underline{\epsilon}} (\uu) : \underline{\underline{\epsilon}} (\bar\uu ) \dr\xx  \geqslant \alpha \|\underline{\underline{\epsilon}}  (\uu)\|_{0,\Omega_{\mathsf{s}}}^2  \ .
   % $$
%Therefore, using the Korn inequality (which is available since $\Omega_{\mathsf{s}}$ is a smooth domain), we infer : For all $\uu\in\bH^1(\Omega_{\mathsf{s}})$
%$$ 
% \int_{\Omega_{\mathsf{s}}}  \underline{\underline{C}} \ \underline{\underline{\epsilon}} (\uu) : \underline{\underline{\epsilon}} (\bar\uu )
%  \dr\xx  \geqslant \alpha C \| \uu\|^2_{1,\Omega_{\mathsf{s}}} -\alpha c  \| \uu\|^2_{0, \Omega_{\mathsf{s}}}\ .
%$$
%\end{rem}
Combining the previous inequality \eqref{korn} and the trace inequality
\begin{equation}
\label{TI}
\forall \uu\in \bH^1(\Omega_{\mathsf{s}})\ ,\quad\|\uu\cdot \nn\|_{0,\Gamma}\leqslant C_{1} \|\uu \|_{1,\Omega_{\mathsf{s}}}
\end{equation}
we infer for $m$ large enough
$$
  \alpha C \|\uu_m\|^2_{1,\Omega_{\mathsf{s}}}  
-\eps_{m} \beta \omega^2\rho_{\sf}  \| \uu_{m} \cdot \nn \|^2_{0,\Gamma}
+\frac{\eps_{m}^3}{3}\omega^2\rho_{\sf} \|  \nabla_{\Gamma}( \uu_{m} \cdot \nn) \|_{0,\Gamma}^2 
   \leqslant C_{2} +  \| \ff_m \|_{0,\Omega_{\mathsf{s}}} \, ,
    $$
 with a constant $\beta >0$ such that $ \cJ_{\eps}\leqslant\beta$ for $\eps$ small enough (remind $\cJ_{\eps}\to 1$ as $\eps \to 0$), and $C_{2}=\alpha c +\omega^2 \rho $. Using again the trace inequality \eqref{TI}, we obtain
 $$
  (\alpha C -\eps_{m}\beta\omega^2\rho_{\sf}C_{1}^2 )  \|\uu_m\|^2_{1,\Omega_{\mathsf{s}}}  
  +\frac{\eps_{m}^3}{3}\omega^2\rho_{\sf} \|  \nabla_{\Gamma}( \uu_{m} \cdot \nn) \|_{0,\Gamma}^2 
    \leqslant C_{2} + \| \ff_m \|_{0,\Omega_{\mathsf{s}}} \,.
  $$
Then, using \eqref{5E4a}, there holds: $\alpha C -\eps_{m}\beta\omega^2\rho_{\sf}C_{1}^2 >0$, for $m$ large enough. According to \eqref{5E4c} we infer that the sequence $\{\uu_m\}$, resp. $\{(\eps_{m})^{\frac32}\nabla_{\Gamma} (\uu_m\cdot\nn)\}$, is bounded in $\bH^1(\Omega_{\mathsf{s}})$, resp. in  $\bL^2(\Gamma)$:
\begin{subequations}
\begin{gather}
\label{5E12}
   \|\uu_m\|_{1,\Omega_{\mathsf{s}}}  
   \leqslant C \,, \\
\label{5E13'}
  (\eps_{m})^{\frac32}\|  \nabla_{\Gamma}( \uu_{m} \cdot \nn) \|_{0,\Gamma}   \leqslant C \ . 
\end{gather}
\end{subequations}

Another consequence of  \eqref{5E12} is that the sequence $\{\uu_m \cdot\nn\}$ is bounded in $\H^{\frac12}(\Gamma)$. 
\begin{equation}
\label{5E13}
   \|\uu_m\cdot\nn\|_{\frac12,\Gamma}  
   \leqslant C \ .
\end{equation}

%According to \eqref{TI}-\eqref{5E12} we infer that the sequence $\{\uu_m\cdot\nn\}$ is bounded in $\L^2(\Gamma)$ :
%\begin{equation}
%\label{5E13}
 %  \|\uu_m\cdot\nn\|_{0,\Gamma}  
  % \leqslant C \ ,
%\end{equation}

\subsubsection{Limit of the sequence and conclusion}

The domain $\Omega_{\mathsf{s}}$ being bounded, the embedding of $\bH^{1}(\Omega_{\mathsf{s}})$ in $\bL^2(\Omega_{\mathsf{s}})$ is compact. Hence as a consequence of \eqref{5E12}, using the Rellich Lemma we can extract a subsequence  of $\{\uu_m\}$ (still denoted by $\{\uu_m\}$) which is converging in $\bL^2(\Omega_{\mathsf{s}})$, and we can assume that the sequence $\{ \underline{\underline{\nabla}} \uu_m\}$ is strongly converging in $\LL^2(\Omega_{\mathsf{s}})$. As a consequence of  \eqref{5E13}, up to the extraction of a subsequence, we can assume that the sequence $\{\uu_m \cdot\nn\}$ is strongly converging in $\L^2(\Gamma)$
: We deduce that there is $\uu\in \bL^2(\Omega_{\mathsf{s}})$ such that
\begin{equation}
\left\{
   \begin{array}{lll}
   \label{4E6}
   \underline{\underline{\epsilon}} (\uu_m) \rightharpoonup  \underline{\underline{\epsilon}} (\uu )  \quad &\mbox{in} \quad \LL^2(\Omega_{\mathsf{s}})
\\
   \uu_m \rightarrow \uu   \quad &\mbox{in} \quad \bL^2(\Omega_{\mathsf{s}})
   \\
    \uu_m \cdot\nn\rightarrow \uu\cdot\nn   \quad &\mbox{in} \quad \L^2(\Gamma) \ .
   \end{array}
\right.
\end{equation}
A consequence of the strong convergence in $\bL^2(\Omega_{\mathsf{s}})$ and \eqref{5E4b} is that $\|\uu\|_{0,\Omega_{\mathsf{s}}}=1$. As a consequence of \eqref{5E13'},  we can extract a subsequence  of 
$\{(\eps_{m})^{\frac32}\nabla_{\Gamma} (\uu_m\cdot\nn)\}$ 
(still denoted by $\{(\eps_{m})^{\frac32}\nabla_{\Gamma} (\uu_m\cdot\nn)\}$) which is weakly converging to a function $\ttt\in \bL^2(\Gamma)$ 
\begin{equation}
\label{4E6'}
(\eps_{m})^{\frac32}\nabla_{\Gamma} (\uu_m\cdot\nn) \rightharpoonup  \ttt \quad \mbox{in} \quad \bL^2(\Gamma)\ .
\end{equation}

Using Assumption \ref{SA}, we are going  to prove that $\uu=0$, which will contradict $\|\uu\|_{0,\Omega_{\mathsf{s}}}=1$, and finally prove estimate \eqref{UE2}. Let $\vv \in \bV^3$ be a test function in \eqref{EVm} 
\begin{multline}
\label{EVm'}
  \int_{\Omega_{\mathsf{s}}}\left(  \underline{\underline{\sigma}}  (\uu_{m}) : \underline{\underline{\epsilon}}  (\bar{\vv})  -\omega^2 \rho \uu_{m}  \bar{\vv}  \right)\;\dr\xx  \
   -\eps_{m}\omega^2\rho_{\sf} \int_{\Gamma} \cJ_{\eps_m}   \uu_{m} \cdot \nn\ \bar{\vv} \cdot \nn \,\dr \sigma 
   \\
    +\frac{\eps_{m}^3}{3}\omega^2\rho_{\sf} \int_{\Gamma}  \nabla_{\Gamma}( \uu_{m} \cdot \nn) \nabla_{\Gamma}(\bar{\vv} \cdot \nn) \,\dr \sigma 
 =   -  \int_{\Omega_{\mathsf{s}}} \ff_{m} \cdot  \bar{\vv} \;\dr\xx  \   .
\end{multline}
According to \eqref{4E6}, \eqref{4E6'}, and \eqref{5E4a}, taking limits as $m\rightarrow +\infty$, there holds 
%the two integrals over $\Gamma$ on the left hand side tends to $0$
$$
  \eps_{m} \int_{\Gamma} \cJ_{\eps_m}   \uu_{m} \cdot \nn\ \bar{\vv} \cdot \nn \,\dr \sigma  \rightarrow 0  
  \quad \mbox{and} \quad
  \frac{\eps_{m}^3}{3}\int_{\Gamma}  \nabla_{\Gamma}( \uu_{m} \cdot \nn) \nabla_{\Gamma}(\bar{\vv} \cdot \nn) \,\dr \sigma \rightarrow 0 \ .
$$
Hence, according to \eqref{4E6}, \eqref{5E4a} and \eqref{5E4c}, taking limits as $m\rightarrow +\infty$, we deduce from the previous equalities \eqref{EVm'} $\uu \in \bH^1(\Omega_{\mathsf{s}})$ satisfies for all $\vv \in \bH^1(\Omega_{\mathsf{s}})$:
\begin{equation*}
%\label{E20}
  \int_{\Omega_{\mathsf{s}}}\left(   \underline{\underline{\sigma}} (\uu) : \underline{\underline{\epsilon}}  (\bar{\vv})  -\omega^2 \rho \uu \bar{\vv}  \right)\;\dr\xx 
  % \ -\eps_{m}\omega^2\rho_{\sf} \int_{\Gamma}    \uu_{m} \cdot \nn\ \bar{\vv} \cdot \nn \,\dr s 
= 0\ .
%\int_{\Gamma}  \hh_{m} \ \bar{\vv} \cdot\nn \,\dr s \   .
\end{equation*}
Integrating by parts we find that $\uu$ satisfies the problem 
\begin{equation*}
 \left\{
   \begin{array}{lll}
\nabla \cdot  \underline{\underline{\sigma}}( \uu)+\omega^2 \rho \uu =0
 \quad&\mbox{in}\quad \Omega_{\mathsf{s}}
\\[0.5ex]
\TT(\uu) =0 \quad &\mbox{on} \quad \Gamma \ .
   \end{array}
    \right.
\end{equation*}
 %Hence $\uu\in\bH^1(\Omega_{\mathsf{s}})$ is solution of problem \eqref{H}. 
 By Assumption \ref{SA}, we deduce
\begin{equation*}
  \uu=0 \quad \mbox{in} \quad \Omega_{\mathsf{s}}\ ,
\end{equation*}
which contradicts $\|\uu\|_{0,\Omega_{\mathsf{s}}}=1$ and ends the proof of Lemma \ref{L1}.

\subsection{Proof of error estimates}
\label{SEe}
In this section we prove the Theorem \ref{CVec}. Since the problem \eqref{Ebck} is of Fredholm type, it is sufficient to show that any solution $\uu_\eps^{k}$ of  \eqref{Ebck} satisfies the error estimate \eqref{UEec}
\begin{equation*}
\label{errork}
\|  \uu_{\eps}- \uu_\eps^k \|_{1,\Omega_{\mathsf{s}}} \leqslant C \eps^{k+1} \ .
\end{equation*}

We prove hereafter the estimate \eqref{UEec} in two steps, Sec. \ref{Si} and Sec. \ref{Sii}.  % \iti1 and \iti2.

\subsubsection{Step A}
\label{Si} The first step consists to derive an expansion of $\uu_\eps^k$ and to show that the truncated expansions of $\uu_\eps^k$  and $\uu_{\eps}$ coincide up to the order $\eps^k$:
\begin{gather}
%\label{AEu}
   \uu_{\eps} = \uu_0 + \eps\uu_1 + \eps^2 \uu_2+ \cdots+\eps^k \uu_{k} +  \mathbf{ r}_{\eps}^{k}   %\mathcal{O}(\eps^3)
   \, ,\\
   \label{AEuk}
  \uu_{\eps}^k = \uu_0 + \eps\uu_1 + \eps^2 \uu_2+ \cdots+\eps^k \uu_{k} +  \tilde{\mathbf{ r}}_{\eps}^{k}\ .
\end{gather}

Hereafter, we justify the expansion \eqref{AEuk}. By construction, $\uu_{\eps}^k$ admits an expansion
$$
 \uu_{\eps}^k = \vv_0 + \eps\vv_1 + \eps^2 \vv_2+ \cdots+\eps^k \vv_{k} +  \tilde{\mathbf{ r}}_{\eps}^{k}\ 
$$
where each term $\vv_n$, for $0\leqslant n\leqslant k$, satisfies the problem \eqref{EAn-} as well as  the term $\uu_n$. Using the spectral Assumption \ref{SA},  we infer that  for all $0\leqslant l\leqslant k$,  $\vv_{n}=\uu_{n}$ in $\Omega_{\s}$, and the expansion \eqref{AEuk} holds.

%The expansion \eqref{AEuk} is obtained as a consequence of the definition of $\uu_{\eps}^k$ together with  the Assumption \ref{SA}. Precisely, using this assumption we can use a uniqueness argument  to ensure that each term in the expansions of $  \uu_{\eps}$ and $  \uu^k_{\eps}$ coincide up to the order $k$.  
Hence, 
\begin{equation*}
%\label{e}
\|  \uu_{\eps}- \uu_\eps^k \|_{1,\Omega_{\mathsf{s}}}
=\|  \mathbf{ r}_{\eps}^{k} -\tilde{\mathbf{ r}}_{\eps}^{k} \|_{1,\Omega_{\mathsf{s}}}
%  \leqslant C \eps^{k+1} 
\ .
\end{equation*}
The estimate of the remainder $\mathbf{ r}_{\eps}^{k} $ is already proved in Thm \ref{TErem} (Sec. \ref{SEER}) : $\|  \mathbf{ r}_{\eps}^{k} \|_{1,\Omega_{\mathsf{s}}}\leqslant C \eps^{k+1}$. In the next  step, we prove estimates for the remainder $\tilde{\mathbf{r}}_{\eps}^{k} $.
 
\subsubsection{Step B}
\label{Sii}
 According to \eqref{AEuk},  the remainder $\tilde{\mathbf{r}}_{\eps}^{k}$ satisfies the elastic equation in $\Omega_{\mathsf{s}}$. We apply the  operator $\TT+\BB_{k,\eps}$ (where $\BB_{k,\eps} (\uu):=\B_{k,\eps}(\uu\cdot\nn) \nn $ )  to the remainder $\tilde{\mathbf{r}}_{\eps}^{k}$. Then, we prove hereafter that
$$
\TT( \tilde{\mathbf{ r}}_{\eps}^{k}) %-\eps \omega^2 \rho_{\sf}  
+\B_{k,\eps}( \tilde{\mathbf{ r}}_{\eps}^{k}\cdot\nn)\nn 
%\cdot \nn \, \nn
= \mathcal{O}(\eps^{k+1})
%  \hh_{k,\eps}
 \quad \mbox{on} \quad \Gamma\ . 
$$
Since $ \uu_{\eps}^0 = \uu_0$, then $\tilde{\mathbf{ r}}_{\eps}^{0}=0$%(i.e. the expansion \eqref{AEuk} is exact at the order $0$)
. Relying on the construction of equivalent conditions detailed in section \ref{CECs},  there holds
 %for $k=1,2$
 %$$
 %\TT( \tilde{\mathbf{ r}}_{\eps}^{k+1}) -\eps \omega^2 \rho_{\sf} \tilde{\mathbf{ r}}_{\eps}^{k+1}\cdot \nn \, \nn= \mathcal{O}(\eps^{k+1}) :
 %$$
 \begin{align*}
%&\TT( \tilde{\mathbf{ r}}_{\eps}^{1})&= &\quad 0 \quad \mbox{on} \quad \Gamma  
%\\
&\TT( \tilde{\mathbf{ r}}_{\eps}^{1}) +\B_{1,\eps}( \tilde{\mathbf{ r}}_{\eps}^{1}\cdot\nn)\nn =\eps^2 \omega^2 \rho_{\sf}\uu_{1}\cdot \nn \,\nn% \quad \mbox{on} \quad \Gamma 
\\[0.8ex]
&\TT( \tilde{\mathbf{ r}}_{\eps}^{2}) +\B_{2,\eps}( \tilde{\mathbf{ r}}_{\eps}^{2}\cdot\nn)\nn 
%-\eps \omega^2 \rho_{\sf}\left(1 +\eps \cH \right)\tilde{\mathbf{ r}}_{\eps}^{2}\cdot \nn \, \nn 
%-\eps \omega^2 \rho_{\sf} \tilde{\mathbf{ r}}_{\eps}^{3}\cdot \nn \, \nn
= \eps^3 \omega^2 \rho_{\sf}\left(\left(1 -\eps \cH \right)\uu_{2} +\cH \uu_{1}\right)\cdot \nn\,\nn
%\quad \mbox{on} \quad \Gamma\ . 
\\[0.8ex]
&\TT( \tilde{\mathbf{ r}}_{\eps}^{3})  +\B_{3,\eps}( \tilde{\mathbf{ r}}_{\eps}^{3}\cdot\nn)\nn  =
\\[0.8ex]
& \eps^4 \omega^2 \rho_{\sf}  
\left[\left( \uu_{3}-  \cH  (\uu_{2}+ \eps \uu_{3})\right)\cdot \nn+\frac13\left( 4 \cH^2 -\cK  +\Delta_{\Gamma} +\kappa^2 \Id \right)
\left( (\uu_{1}+ \eps \uu_{2}+\eps^2 \uu_{3})\cdot \nn\right)
\right] \nn  
\end{align*}
on $\Gamma$. According to the estimate \eqref{UE2}, we infer the uniform estimate 
\begin{equation*}
\label{UErk+1}
\|  \tilde{\mathbf{ r}}_{\eps}^{k} \|_{1,\Omega_{\mathsf{s}}} \leqslant C \eps^{k+1}
\  ,
\end{equation*}
which ends the proof of Theorem \ref{CVec}.

\section*{Acknowledgements}

The author thanks very warmly Monique Dauge for her well-considered suggestions.

%%-----------------------------
%%      your bibliography
%%-----------------------------
\bibliographystyle{plain}
\bibliography{biblio}

\end{document}